\documentclass{svmult}
\usepackage{amsmath}
\usepackage{amssymb}
\usepackage{amsfonts,amsmath}
\usepackage{latexsym}
\usepackage{euscript}
\usepackage{graphicx}

\def\R{\rz R}

\def\E{\rz E}

\def\shc{{\cal C}}

\def\shf{{\cal F}}
\def\shg{{\cal G}}

\def\shm{{\cal M}}

\def\shm{{\cal M}}
\def\shl{{\cal L}}

\def\shy{{\cal Y}}

\newcommand{\rz}[1]{\mathbb{#1}}

\newcommand{\dis}{\displaystyle}

\newcommand{\eps}{\varepsilon}
\newcommand{\beqnar}{\begin{eqnarray*}}
\newcommand{\eeqnar}{\end{eqnarray*}}
\newcommand{\ba}{\begin{array}}
\newcommand{\ea}{\end{array}}

\newcommand{\halb}{\frac{1}{2}}

\newcommand{\noi}{\noindent}

\begin{document}

\title*{ELEMENTS OF STOCHASTIC CALCULUS VIA REGULARISATION}

\toctitle{Elements of stochastic calculus via regularisation}

\titlerunning{Stochastic calculus via regularisation}

\author{Francesco Russo (1) and Pierre Vallois (2)}

\tocauthor{Francesco Russo and Pierre Vallois}

\authorrunning{Francesco Russo and Pierre Vallois}

\institute{(1) Universit\'{e} Paris 13\\
           Institut Galil\'{e}e, Math\'ematiques\\
          99 avenue J.B.~Cl\'{e}ment\\
           F-93430 Villetaneuse, France\\
           \texttt{e-mail: russo@math.univ-paris13.fr}\\
            (2) Universit\'e Henri Poincar\'e\\
          Institut de Math\'ematiques Elie Cartan\\
           B.P. 239\\
       F-54506 Vand\oe uvre-l\`es-Nancy Cedex, France\\
       \texttt{e-mail: vallois@iecn.u-nancy.fr}}

\maketitle

\begin{abstract}

This paper first summarizes the foundations of stochastic calculus
via regularization and constructs
through this procedure It\^o and Stratonovich integrals.
In the second part,  a survey and  new results are presented
in relation with finite quadratic variation processes,  Dirichlet and
weak Dirichlet processes.

\

\noi{\bf Keywords:} Integration via regularization, weak Dirichlet
processes, covariation, It\^o formulae.

\

\noi{\bf MSC 2000:} 60H05, 60G44, 60G48

\end{abstract}

\section{Introduction}

Stochastic integration via regularization is a technique of
integration developed in a series of papers by the authors
starting from \cite{rv0}, continued in \cite{rv1, rv93, rv2, rv4,
rv96} and later carried out by other authors, among them
\cite{rvw, er, er2, wo1, wo2, wo3, za,
flru1, flru2, frw1, frw2,
gr}. Among some recent applications to finance,
we refer for instance to \cite{lnn, oks}.

This approach constitutes a counterpart of a discretization
approach initiated by F\"ollmer (\cite{fo})  and continued by
many authors, see for instance \cite{ber, fps, fp, erv, eisen,
gp}.

The two theories run parallel and, at the axiomatic level,  
almost all  the results we
obtained via regularization can  essentially be translated in the
language of discretization.

The advantage of using regularization lies in the fact that this
approach is natural and relatively  simple, and easily connects to other
approaches. We now list some typical features of stochastic calculus via
regularization.

\begin{itemize}
\item Two fundamental notions are the 
quadratic variation of a process, see
Definition \ref{D23} and the forward integral, see Definition \ref{D20}.
Calculus via regularization is first of all a calculus related to finite
quadratic variation processes, see section~\ref{S6}.
 It\^o integrals with respect to  continuous semimartingales
can be defined through  forward integrals, see Section~\ref{S5};
this makes classical stochastic calculus  appear as a 
particular instance of  calculus via regularization.
Let  the integrator be a classical Brownian motion $W$ and the
integrand a measurable adapted process $H$ such that 
$\int_0^T H^2_t dt < \infty$ a.s.,
where a.s. means almost surely.
We will show in section~\ref{SS54B} that the forward integral 
$\int_0^\cdot H d^-W$ coincides with the It\^o integral
$\int_0^\cdot H dW$.
On the other hand, the  discretization approach constitues a sort
of Riemann-Stieltjes type integral and  only allows  
integration of processes that are not too
irregular, see  Remark \ref{RITOV2}.

 \item Calculus via regularization constitutes a bridge between
non causal and causal calculus operating through substitution
formulae, see subsection~\ref{SSubsti}. A precise link between
forward integration and  the theory of enlargement
of filtrations may be given, see \cite{rv1}.
 Our integrals can be  connected to the well-known
Skorohod type integrals, see again \cite{rv1}.

\item With the help of symmetric integrals a calculus with respect
to processes with a variation higher than 2 may be developed.
For instance fractional Brownian motion is the prototype of
such processes.

\item This stochastic calculus constitutes somehow a 
barrier separating the pure pathwise calculus in the
sense of T. Lyons and coauthors, see e.g. \cite{lyons, lq, lejay, gubi},
and any stochastic calculus  taking into account an underlying
probability,  see   Section~\ref{Slast}.
\end{itemize}

This paper will essentially focuse on  the first
item.

The paper is organized as follows. First, in Section
\ref{Stochint}, we recall the basic definitions  and properties of
forward, backward, symmetric integrals and covariations. Justifying
the related definitions and properties needs no particular
effort.
A significant example is the Young integral, see \cite{young}.
In Section~\ref{S5} we redefine It\^o integrals in the spirit of
integrals via regularization and we prove
some typical properties. We  essentially define
It\^o integrals as forward integrals in a subclass
and we then extend this definitionthrough functional analysis methods.
 Section~\ref{S6}
 is devoted to finite quadratic variation
processes. In particular we establish  $C^1$-stability properties
and an It\^o formula of $C^2$-type. Section~\ref{S8}
provides  some survey material with  new results 
related to the  class of weak Dirichlet
processes introduced by \cite{er} with later developments
discussed by \cite{gr, cjms}. Considerations about  
It\^o formulae under $C^1$-conditions
are discussed as well.

\section{Stochastic integration via regularization} \label{Stochint}

\subsection{Definitions and fundamental properties}\label{Stochint1}

In this paper $T$ will be a fixed positive real number.
By convention, any real continuous function $f$ 
defined either on $[0,T]$
or $\R_+$ will be prolongated (with the 
same name) to the real line, setting

\begin{equation}\label{E5PV1}
    f(t) = \left \{
\begin{array} {ccc}
 f(0) &{\rm if} & t\le 0 \\
f(T) &{\rm if} & t > T.
\end{array}
\right.
\end{equation}

Let ${(X_t)}_{t \geq 0 }$ be a continuous process and ${(Y_t)}_{ t
\geq 0}$ be a process with paths in $L^1_{loc}(\mathbb{R}_+)$,
i.e. for any $a>0$,  $ \displaystyle \int_0^a \vert Y_t \vert dt
< \infty $ a.s.

Our generalized stochastic integrals and covariations will be
defined through a regularization procedure. More precisely, let
$I^-(\varepsilon , Y, dX)$ (resp. $I^+(\varepsilon , Y, dX)$,
 $I^0(\varepsilon , Y, dX)$ and $C(\varepsilon , Y, X)$) be the
$\varepsilon$-forward integral (resp. $\varepsilon$-backward
integral, $\varepsilon$-symmetric integral and
$\varepsilon$-covariation):
\begin{equation}\label{SIVR1n}
I^-(\varepsilon , Y, dX)(t)=\int _0^tY(s)
\frac{X(s+\varepsilon)-X(s)}{\varepsilon}ds; \quad  t \geq 0,
\end{equation}
\begin{equation}\label{SIVR1An}
I^+(\varepsilon , Y, dX)(t)=\int _0^tY(s)
\frac{X(s)-X(s-\varepsilon}{\varepsilon}ds; \quad  t \geq 0,
\end{equation}
\begin{equation}\label{SIVR1Bn}
I^0(\varepsilon , Y, dX)(t)=\int _0^tY(s)
\frac{X(s+\varepsilon)-X(s-\varepsilon )}{2\varepsilon}ds;
\quad  t \geq 0,
\end{equation}
\begin{equation}\label{SIVR1Cn}
C(\varepsilon , X, Y)(t)=\int _0^t
\frac{\big(X(s+\varepsilon)-X(s)\big)\big(Y(s+\varepsilon)
-Y(s)\big)}{\varepsilon}ds; \quad  t \geq 0 .
\end{equation}
Observe that these four processes are continuous.

\begin{definition} \label{D20}\begin{description}
    \item {1)} A family of processes $(H^{(\varepsilon)}_t)_{
    t \in  [0,T]}$ is said to converge  to $(H_t)_{
    t \in [0,T] }$ in the {\bf  ucp sense},
if  $\displaystyle \sup_{0\leq
    t\leq T}|H^{(\varepsilon)}_t-H_t|$ goes to $0$ in probability,
    as $\varepsilon \rightarrow 0$.

    \item{2)}
    Provided the corresponding limits exist in the ucp
 sense, we define the following integrals and covariations by the
following formulae

\begin{description}

\item{a)} {\rm \bf Forward integral:}  $\displaystyle \int_0^t Y
d^- X = \lim_{\varepsilon \rightarrow 0^+} I^-(\varepsilon , Y,
dX)(t)$.

\item{b)} {\rm \bf Backward integral:}  $ \displaystyle \int_0^t Y
d^+ X = \lim_{\varepsilon \rightarrow 0^+} I^+(\varepsilon , Y,
dX)(t).$

\item{c)} {\rm \bf Symmetric integral:}   $\displaystyle \int_0^t
Y d^\circ X = \lim_{\varepsilon \rightarrow 0+} I^\circ (\varepsilon ,
Y, dX)(t)$.

\item{d)} {\rm \bf Covariation:}  $\displaystyle[X,Y]_t =
\lim_{\varepsilon \rightarrow 0^+} C(\varepsilon , X, Y)(t)$.
 When $X = Y$ we often put $[X] = [X,X]$.

\end{description}

\end{description}

\end{definition}



\begin{remark} \label{R21} Let $X, X', Y, Y'$ be four
processes with $X,X'$ continuous and $Y, Y'$ having paths in
$L^1_{loc}(\mathbb{R}_+)$.   $\star$ will stand for one of the three
symbols $-$, $+$ or $\circ$.
\begin{enumerate}
\item $ (X,Y) \mapsto \int_0^\cdot Y d^\star X$ and
$ (X,Y) \mapsto [X,Y]$ are bilinear operations.

\item  The covariation of continuous processes is a symmetric
operation.

\item   When it exists, $[X]$ is  an increasing process.

 \item   If $\tau$ is a random time, $
{[X^\tau, X^\tau]}_t = {[X,X]}_{t\wedge \tau} $  and
$$ \int_0^t Y 1_{[0,\tau]} d^\star X = \int_0^t Y d^\star X^\tau =
 \int_0^t Y^\tau d^\star X^\tau =  \int_0^{t \wedge \tau} Y d^\star X,
$$
where  $X^\tau$ is the process $X$
stopped at time $\tau$, defined by $X^\tau_t = X_{t\wedge \tau} $.

\item  If $\xi$ and $\eta$ are two fixed r.v., 
$\displaystyle \int_0^\cdot (\xi Y_s) d^\star (\eta X_s)
 = \xi \eta \int_0^\cdot Y_s d^\star X_s. $

\item  Integrals via regularization also have the following
localization property.
 Suppose
that  $X_t = X'_t , Y_t = Y'_t, \forall t \in [0,T]$ on some
subset $\Omega_0$ of $\Omega.$ Then
$$ 1_{\Omega_0} \int_0^t  Y_s d^\star X_s
 = 1_{\Omega_0} \int_0^t Y'_s d^\star X'_s, \quad t \in [0,T]. $$

 \item  If $Y$ is an elementary process of
the type $\displaystyle Y_t = \sum_{i= 1}^N A_i 1_{I_i}, $ where
$A_i$ are random variables and $(I_i)$  a family of real
intervals with end-points $a_i < b_i$, then
$$ \int_0^t  Y_s d^\star  X_s =
 \sum_{i= 1}^N A_i (X_{b_i \wedge t} - X_{a_i \wedge t}).
$$

\end{enumerate}
 \end{remark}

\begin{definition} \label{D23}
\begin{description}
\item {1)} If $[X]$ exists, $X$ is said to be a {\bf finite quadratic
variation process} and $[X]$ is called the {\bf quadratic variation} of
$X$.

\item {2)} If $[X] = 0$, $X$ is called a {\bf zero quadratic
variation process}.

\item {3)} A vector $(X^1,\ldots,X^n)$ of continuous processes is
said to have all its {\bf mutual covariations }
 if
$[X^i, X^j]$ exists for all $1 \leq i$, $j \leq n$.

We will also use  the terminology {\bf bracket} instead of covariation.
\end{description}
\end{definition}

\begin{remark} \label{R24}
\begin{description}
    \item {1)} If $(X^1,\ldots,X^n)$ has all its mutual covariations, then
\begin{equation} \label{EBil}
[X^i + X^j, X^i + X^j] = [X^i, X^i] + 2 [X^i, X^j] + [X^j, X^j].
\end{equation}
From the previous equality, it follows that
 $[X^i, X^j]$ is the difference of two increasing processes,
having therefore bounded variation; consequently the  bracket is  a classical
integrator in the Lebesgue-Stieltjes sense.

    \item {2)} Relation
(\ref{EBil}) holds as soon as three brackets among  the four
exist. More generally, by convention,  an identity of the type $I_1 + \cdots + I_n
= 0$ has the following meaning: if   $n - 1$ terms among the $I_j$
exist,  the remaining one also makes sense and the identity holds
true.

 \item {3)} We will see later, in Remark \ref{R816}, that there exist
processes $X$ and $Y$  such that $[X,Y]$ exists but does not have
finite variation; in particular $(X,Y)$ does not have all its
mutual brackets.

\end{description}
\end{remark}

The properties below follow elementarily from 
the definition of  integrals via regularization.

\begin{proposition} \label{P24}
Let $X = {(X_t)}_{t \geq 0}$ be a continuous process and $Y =
{(Y_t)}_{t \geq 0} $ be a  process with  paths in
$L^1_{loc}(\mathbb{R}_+)$. Then
\begin{description}
\item{1)} $[X,Y]_t = \displaystyle \int_0^t Yd^+ X - \int_0^t Yd^-
X$.

\item{2)} $\displaystyle \int_0^t Y d^\circ X = \frac{1}{2} \left(
\int_0^t Yd^+ X + \int_0^t Yd^-X\right)$.

 \item{3)}  {\bf Time reversal}.  Set $\hat X_t = X_{T-t}, t \in  [0,T]$. 
 Then
\begin{enumerate}
\item $\displaystyle \int_0^t Yd^{\pm} X = - \int_{T-t}^T \hat Y
d^{\mp} \hat X  , \quad  0 \leq t\leq T$ ;

\item $\displaystyle \int_0^t Yd^\circ X = - \int_{T-t}^T \hat Y
 d^\circ \hat X , \quad  0 \leq t\leq T$;

\item $[\hat X,\hat Y]_t = [X,Y]_T - [X,Y]_{T-t} , \quad  0 \leq
t\leq T$.
\end{enumerate}

\item{4)}  {\bf Integration by parts}. If  $Y$ is continuous,
\begin{eqnarray*}
X_t Y_t &=& X_0 Y_0 + \int_0^t Xd^- Y + \int_0^t Yd^+ X\\
        &=& X_0 Y_0 + \int_0^t Xd^- Y + \int_0^t Yd^- X + [X,Y]_t.
\end{eqnarray*}
\item{5)} {\bf Kunita-Watanabe inequality}. If $X$ and $Y$ are
finite quadratic variation processes,
$$ \bigl| [X,Y] \bigr| \le \big\{[X]\  [Y] \big\}^{1/2} .$$

\item{6)} If $X$ is a finite quadratic variation process and $Y$
is a zero quadratic variation process
 then $(X,Y)$ has all its mutual brackets and $[X,Y] = 0$.

\item{7)} Let $X$ be a bounded variation process and
 $Y$ be a process with locally bounded paths,  and at most countably
 many discontinuities. Then
\begin{description}
\item{a)} $\displaystyle \int_0^t Yd^+ X = \int_0^t Yd^- X =
\int_0^t YdX$, where $\displaystyle  \int_0^t YdX$ is a 
Lebesgue-Stieltjes integral.
 \item{b)} $[X,Y] = 0$. In particular a bounded
variation and continuous process is a zero quadratic variation process.
\end{description}
\item{8)} Let $X$ be an absolutely continuous process and
 $Y$ be a process with locally bounded paths. Then
$$ \displaystyle \int_0^t Yd^+ X = \int_0^t Yd^- X =
\int_0^t Y X' ds.$$
\end{description}
\end{proposition}

\begin{remark}\label{r25bis}
If $Y$ has uncountably many discontinuities,
7) may fail.
Take for instance $Y = 1_{{\rm supp}\, dV}$,
where $V$ is an increasing continuous
function such that $V'(t) = 0$ a.e. (almost everywhere) with respect to Lebesgue measure.
Then
$ Y = 0 $ Lebesgue a.e., and $Y =1, dV$ a.e.
Consequently
$$ \int_0^t Y d V = V(t) - V(0), \quad   I^- (\varepsilon, Y, dV)(t) = 0   \quad    \int_0^t Y d^-V = 0.$$
\end{remark}

\begin{remark} \label{R25}
Point 2) of Proposition \ref{P24} states that the symmetric
integral is the average of the forward and backward integrals.
\end{remark}

\noi\textit{Proof} of Proposition \ref{P24}.
Points 1), 2), 3), 4) follow immediately from the definition. For
illustration, we only prove 3); operating a change of variable $u
= T-s$, we obtain
$$
 \int_0^t Y_s \frac{X_s - X_{s-\varepsilon}
}{\varepsilon} ds
 = - \int_{T-t}^{T} \hat Y_u \frac{ \hat X_{u+\varepsilon}
- \hat X_u}{\varepsilon} du , \quad 0 \leq t \leq T.
$$
Since $X$ is continuous, one can take the limit of both members and
the result follows.

5) follows by Cauchy-Schwarz inequality which says that
\begin{eqnarray*}
&& \frac{1}{\varepsilon} \left | \int_0^t
(X_{s+\varepsilon}-X_s)\,
(Y_{s+\varepsilon}-Y_s) ds\right | \\
&&\quad \leq \left\{ \frac{1}{\varepsilon} \int_0^t
(X_{s+\varepsilon}-X_s)^2 ds  \ \frac{1}{\varepsilon} \int_0^t
(Y_{s+\varepsilon}-Y_s)^2 ds\right\}^\halb .
\end{eqnarray*}

6) is a consequence of 5).

7) Using Fubini, one has
\begin{eqnarray*}
 \frac{1}{\varepsilon}  \int_0^t Y_s(X_{s+\varepsilon}-X_s)ds &=&
\frac{1}{\varepsilon} \int_0^t ds\, Y_s \int_s^{s+\varepsilon} dX_u\\
 &=&\int_0^{t+\varepsilon} dX_u\;\frac{1}{\varepsilon} 
 \int_{u-\varepsilon}^{u \wedge t}
Y_s ds .
\end{eqnarray*}

Since the jumps of $Y$ are at most countable,
 $\displaystyle
\frac{1}{\varepsilon} \int_{u -\varepsilon}^u Y_s ds
\rightarrow Y_u,  d \vert X \vert$ a.e.   
where $\vert X \vert$ denotes the total
variation of $X$.
Since
 $t\rightarrow Y_t$ is locally bounded,
Lebesgue's convergence theorem implies that  $\displaystyle
\int_0^t Yd^- X = \int_0^t YdX$.

 The fact that $\displaystyle
\int_0^t Yd^+ X = \int_0^t YdX$ follows similarly.

b) is a consequence of point 1).

8) can be reached using similar elementary integration properties.

\subsection{ Young integral in a simplified framework} \label{SYoung}

We will consider the  integral defined by Young  (\cite{young}) in
1936, and implemented in the stochastic framework by 
Bertoin, see \cite{ber2}.
Here we will restrict ourselves to the case when integrand and integrator
are H\"older continuous processes.
As a result, that integral will be shown to coincide with the forward
integral, 
but also with backward and symmetric ones.

\begin{definition}\label{DEXVQ3}\begin{enumerate}
    \item Let ${C}^\alpha$ be the set of H\"{o}lder continuous
functions defined on $[0,T]$, with index $\alpha >0$. Recall that
$f : [0,T]\mapsto \mathbb{R}$ belongs to ${C}^\alpha$ if
$$ N_\alpha(f):=\sup_{0 \leq s,t \leq
T}\frac{|f(t)-f(s)|}{|t-s|^\alpha}<\infty .$$
    \item If  $X,Y : [0,T]\mapsto \mathbb{R}$ are two functions  of class
    $C^1$, the Young integral of $Y$ with respect to $X$ on 
    $[a,b]\subset [0,T]$ is
    defined as :
    $$  \int_a^b Y d^{(y)}X:=\int_a^b Y(t) X'(t)dt, \quad 0\leq a
    \leq b \leq T.$$
\end{enumerate}

\end{definition}

To extend the Young integral to H\"{o}lder functions we need some
estimate of $\displaystyle \int_0^T Y d^{(y)}X$ in terms of
the H\"{o}lder
norms of $X$ and $Y$. More precisely, let $X$ and $Y$ be as in
Definition \ref{DEXVQ3} above; then
 in \cite{fp}, it is proved:
    \begin{equation}\label{EXVQ32}
    \Bigl|\int_a^T (Y - Y(a)) d^{(y)}X \Bigr| 
    \leq C_\rho T^{1+\rho} N_\alpha (X)
    N_\beta (Y), \quad 0 \le a \le T,
\end{equation}
where $\alpha, \beta >0$, $\alpha + \beta>1,\ \rho \in ]0,\alpha
+\beta -1[$, and $C_\rho$ is a universal constant.


\begin{proposition} \label{PEXVQ7b}\begin{enumerate}
    \item The map  $\displaystyle (X,Y)\in C^1([0,T])\times
C^1([0,T])\mapsto \int_0^\cdot Y d^{(y)}X$ with values in
${C}^\alpha$, extends to a continuous bilinear
map from ${C}^\alpha \times {C}^\beta$ to ${\cal
C}^\alpha$. The value of this extension at point $(X,Y)\in {\cal
C}^\alpha \times {C}^\beta$ will  still be denoted by $\displaystyle
\int_0^\cdot Y d^{(y)}X$ and called the {\bf Young integral} of
$Y$ with respect to $X$.
    \item Inequality (\ref{EXVQ32}) is still valid for any
    $X\in {C}^\alpha$ and  $Y \in {C}^\beta$.
\end{enumerate}

\end{proposition}

\begin{proof} \ 1.  Let $X,Y$ be of class $C^1([0,T])$ and
$$F(t)=\int_0^t Y d^{(y)}X=\int_0^t Y(s) X'(s)ds, \quad t\in
[0,T].$$
For any $a,b \in [0,T], \ a<b$, we have
$$F(b)-F(a)=\int_a^b\big(Y(t)-Y(a)\big)d^{(y)}X+Y(a)\big(X(b)-X(a)\big).$$
Then (\ref{EXVQ32}) implies
\begin{equation}\label{EXVQ32b}
    |F(b)-F(a)|\leq C_\rho(b-a)^{1+\rho}N_\alpha (X)
    N_\beta (Y)+\sup_{0 \leq t \leq T}|Y(t)|\ N_\alpha
    (X)(b-a)^\alpha;
\end{equation}
consequently $F\in {C}^\alpha$.

Then the map $\displaystyle (X,Y)\in C^1([0,T])\times C^1([0,T])\mapsto
\int_0^\cdot Y d^{(y)}X$, which is bilinear,  extends to a
continuous bilinear map from ${C}^\alpha \times {\cal
C}^\beta$ to ${C}^\alpha$.

2. is a consequence of point 1.
\end{proof}

Before discussing the relation between  Young integrals and integrals
via regularization, here is useful technical result.

\begin{lemma} \label{LEXVQ2} Let $0<\gamma'<\gamma \leq 1,
\varepsilon>0$. With $Z\in{C}^\gamma$ we associate
$$Z_\varepsilon (t)=\frac{1}{\varepsilon}\int_0^t\big(Z(u+\varepsilon) - Z(u)\big)du, \ t \in [0,T].$$
Then $Z_\varepsilon$ converges to $Z$ in ${C}^{\gamma'}$, as
$\varepsilon\rightarrow 0$.
\end{lemma}

\begin{proof}
For any $0\leq t\leq T$,
$$Z_\varepsilon(t)=\frac{1}{\varepsilon}\int_0^t\big(Z(u+\varepsilon)-Z(u)\big)du=
\frac{1}{\varepsilon}\int_t^{t+\varepsilon}Z(u)du -
\frac{1}{\varepsilon}\int_0^{\varepsilon}Z(u)du.$$

  Setting
$\Delta_\varepsilon(t)=Z_\varepsilon(t)-Z(t)$, we get
$$
\begin{array}{ccl}
 \Delta_\varepsilon(t ) -\Delta_\varepsilon(s)&=& \displaystyle
\frac{1}{\varepsilon}\int_t^{t+\varepsilon}Z(u)du-Z
(t)-\frac{1}{\varepsilon}\int_s^{s+\varepsilon}Z(u)du+Z
(s) \\
   & = & \displaystyle \frac{1}{\varepsilon}\int_t^{t+\varepsilon}\big(Z(u)-Z
(t)\big) du
-\frac{1}{\varepsilon}\int_s^{s+\varepsilon}\big(Z(u)-Z
(s)\big) du, \\
\end{array}
$$
where $0 \leq s \leq t \leq T$.

a) Suppose $0\leq s <s+\varepsilon<t$. The above inequality
implies
$$|\Delta_\varepsilon(t ) -\Delta_\varepsilon(s)|\leq
\frac{1}{\varepsilon}\int_t^{t+\varepsilon}\big|Z(u)-Z
(t)\big| du
+\frac{1}{\varepsilon}\int_s^{s+\varepsilon}\big|Z(u)-Z
(s)\big| du.
$$
%

%
Since $Z \in {C}^\gamma$, then
$$\begin{array}{ccl}
  |\Delta_\varepsilon(t ) -\Delta_\varepsilon(s)| & \leq & \displaystyle\frac{N_\gamma
(Z)}{\varepsilon}\Big(\int_t^{t+\varepsilon}(u-t)^\gamma du +
\int_s^{s+\varepsilon}(u-s)^\gamma du \Big)\\
   & \leq  &\displaystyle  \frac{2N_\gamma
(Z)}{\gamma +1}\varepsilon^\gamma. \\
\end{array}
$$
But $\varepsilon <t-s$, consequently
\begin{equation}\label{EXVQ33}
|\Delta_\varepsilon(t ) -\Delta_\varepsilon(s)|  \leq
\frac{2N_\gamma (Z)}{\gamma
+1} \varepsilon^{\gamma-\gamma'}|t-s|^{\gamma'}.
\end{equation}
b) We now investigate the case $0\leq s < t < s + \varepsilon$. The
difference $\Delta_\varepsilon(t) -\Delta_\varepsilon(s)$ may
be decomposed as follows :
$$\begin{array}{ccl}
  \Delta_\varepsilon(t ) -\Delta_\varepsilon(s) & = &\displaystyle
  \frac{1}{\varepsilon}\int_{s+\varepsilon}^{t+\varepsilon}
  \big(Z(u)-Z(s+\varepsilon)\big) du-
  \frac{1}{\varepsilon}\int_{s}^{t}
  \big(Z(u)-Z(s)\big) du \\
  &&\\
   &  &\displaystyle +\frac{t-s}{\varepsilon}\big(Z(s+\varepsilon)-Z(s)
   \big)+Z(s) - Z(t).\\
\end{array}
$$
Proceeding as in the previous step and using the inequality
$0<t-s<\varepsilon$, we obtain
$$
\begin{array}{ccl}
  |\Delta_\varepsilon(t ) -\Delta_\varepsilon(s)| & \leq &\displaystyle N_\gamma (Z)\Big(
  \frac{2}{\gamma +1}\frac{(t-s)^{\gamma+1}}{\varepsilon}
  +\frac{t-s}{\varepsilon^{1-\gamma}}+(t-s)^\gamma \Big)\\
   & \leq  &\displaystyle 2N_\gamma (Z)\frac{\gamma+2}{\gamma+1} \varepsilon^{\gamma-\gamma'}
   |t-s|^{\gamma'}. \\
\end{array}
$$
At this point, the above inequality and (\ref{EXVQ33})
directly imply that $N_{\gamma'}(Z_\varepsilon -Z)\leq C
\varepsilon^{\gamma-\gamma'}$ and the  claim is finally established.
\end{proof}

In the sequel of this section $X$ and $Y$ will denote stochastic processes.

\begin{remark} \label{EEXVQ1} If $X$ and $Y$ have a.s. H\"older
continuous paths respectively of order $\alpha$ and $\beta$
with $\alpha> 0$, $\beta >0$ and $\alpha + \beta > 1$. 
Then one can easily prove that
$[X,Y]=0$.
\end{remark}

\begin{proposition} \label{PEXVQ7}  Let $ X,Y $ be two real processes
indexed by $[0,T]$ whose paths are respectively a.s. in
 $ {C}^\alpha$ and  ${C}^\beta$, with $\alpha
>0, \beta >0$ and $\alpha + \beta >1$.
Then the three integrals
    $\displaystyle \int_0^\cdot Yd^+ X$,
    $\displaystyle \int_0^\cdot Yd^- X$
    and $\displaystyle \int_0^\cdot Yd^\circ X$ exist
    and coincide with the
    Young integral  $\displaystyle \int_0^\cdot Yd^{(y)} X$.
\end{proposition}

\begin{proof}

 We  establish
that the forward integral coincides with the Young integral.
The equality concerning the two other integrals is a consequence
of Proposition \ref{P24} 1., 2. and Remark \ref{EEXVQ1}.

By additivity we  can suppose, without lost  generality, that
$Y(0)=0$.

Set
$$\Delta _\varepsilon(t):=\int_0
^tYd^{(y)}X-\int_0^tYdX_\varepsilon,$$
where
$$
X_\varepsilon (t)=\frac{1}{\varepsilon}\int_0^t
\big(X(u+\varepsilon)-X(u)\big)du, \ t \in [0,T].
$$
Since $t \mapsto X_\varepsilon (t)$ is of class $C^1([0,T])$, then
$\displaystyle \int_0^tYdX_\varepsilon$ is equal to the Young integral
$\displaystyle \int_0^tYd^{(y)}X_\varepsilon$ and therefore
$$
\Delta _\varepsilon(t)=\int_0
^tYd^{(y)}\big(X-X_\varepsilon\big).
$$
 Let $\alpha '$ be such that : $0<\alpha'<\alpha$ and
 $\alpha'+\beta>1$. Applying inequality (\ref{EXVQ32}) we obtain
$$
\sup_{0\leq t \leq T}|\Delta _\varepsilon(t)| \leq C_\rho
T^{1+\rho}N_{\alpha'} (X-X_\varepsilon)N_{\beta}(Y),\quad \rho \in
]0,\alpha' {+}\beta{-}1[.
$$
Lemma \ref{LEXVQ2} with $Z=X$ and $\gamma=\alpha$ directly implies
that $\Delta _\varepsilon(t)$ goes to $0$, uniformly a.s. on $[0,T]$,
as $\varepsilon\rightarrow 0$,   concluding  the proof of the Proposition.

\end{proof}

\section{ It\^o integrals and related topics} \label{S5}

 The section presents the construction of It\^o integrals
with respect to continuous local martingales; it is  based on
  McKean's idea (see section 2.1 of \cite{mk}),
which fits the spirit of calculus via regularization.

\subsection{Some reminders on martingales theory} \label{S5mart}

In this subsection, we recall basic notions related to martingale theory,
essentially without proofs, except when they help the reader.
For detailed complements, see \cite{ks}, chap. 1.,
in particular for  definition of adapted and progressively measurable
processes.

Let $(\shf_t)_{t \geq 0}$ be a filtration on
the probability space $(\Omega, \shf, P)$
satisfying  the usual conditions, see  Definition 2.25, chap. 1 in \cite{ks}.

An adapted process $(M_t)$ of integrable random variables,
 i.e. verifying $E( |M_t|) < \infty, \ \forall t \geq 0$ is:
\begin{itemize}
\item an $(\shf_t)$-martingale if
 $E( M_t| \shf_s) = M_s, \quad \forall t \geq s$;

\item a $(\shf_t)$- submartingale if
 $E( M_t| \shf_s) \geq M_s, \quad \forall t \geq s$
\end{itemize}

In this paper,  all submartingales (and therefore all
martingales) will be supposed to be continuous.

\begin{remark} \label{R341}
 It follows from the definition that if $(M_t)_{t \geq 0}$
 is a martingale, then $E(M_t) = E(M_0), \ \forall t \geq 0$.
 If $(M_t)_{t \geq 0}$ is a supermartingale (resp. submartingale)
 then $t \longrightarrow E(M_t)$ is decreasing (resp. increasing).
\end{remark}

\begin{definition} \label{D341}
A process $X$ is said to be {\bf square integrable} if
 $E(X_t^2) < \infty$ for each $t \ge 0.$
\end{definition}

When we speak of a martingale without specifying the
$\sigma$-fields, we refer to the {\it canonical} filtration
generated by the process and satisfying the usual conditions.

\begin{definition} \label{D326}
\begin{enumerate}
\item A (continuous) process ${(X_t)}_{t \ge 0}$, is
called a $(\shf_t)$-{\bf local martingale} 
(resp. $(\shf_t)$-{\bf local submartingale})
if there
exists an increasing sequence $(\tau_n)$ of stopping times
such that
 $X^{\tau_n} 1_{\tau_n > 0}  $ is an  $(\shf_t)$-martingale
(resp.  submartingale)
and $ \displaystyle \lim_{n \rightarrow \infty} \tau_n = +\infty$ a.s.
\end{enumerate}
\end{definition}

\begin{remark} \label{R327}

\begin{itemize}
\item An   $(\shf_t)$-martingale  is an 
 $(\shf_t)$-local martingale. A bounded  
  $(\shf_t)$-local martingale is an $(\shf_t)$-martingale.
\item The set of $(\shf_t)$-local martingales  
 is a linear space.
\item If $M$ is an  $(\shf_t)$-local martingale
and $\tau$ a stopping time, then
$M^\tau$ is again an $(\shf_t)$-local martingale.
\item If $M_0$ is bounded, in the definition of a local martingale
one can choose a localizing sequence $(\tau_n)$
such that each $M^{\tau_n}$ is bounded.
\item A convex function of an $(\shf_t)$-local submartingale is an $(\shf_t)$-local submartingale.

\end{itemize}
\end{remark}

\begin{definition} \label{D327}
A process $S$ is called a (continuous)  
$(\shf_t)$-{\bf semimartingale} if it is
the sum of an  $(\shf_t)$-local martingale and an $(\shf_t)$-adapted
continuous bounded variation process.
\end{definition}

A basic decomposition in stochastic analysis is the following.

\begin{theorem} \label{P331}
({\bf Doob decomposition of a submartingale})

Let $X$ be a   $(\shf_t)$-local submartingale.
Then, there is an $(\shf_t)$-local martingale $M$
and an adapted, continuous, and finite variation process $V$ (such
that $V_0 = 0$) with $X = M + V$. The decomposition is unique.
\end{theorem}

\begin{definition} \label{D328}
Let $M$ be an  $(\shf_t)$-local martingale.
We denote by $<M>$ the bounded variation process
featuring in the  Doob decomposition of 
the local submartingale $M^2$.
In particular $M^2 - <M>$ is an   $(\shf_t)$-local martingale.
\end{definition}
In Corollary \ref{CDrob}, we will prove  that $<M>$ coincides with $[M,M]$, so that the skew bracket $<M>$ does not depend on
the underlying filtration.

The following result will be needed in section~\ref{SITOInt}.

\begin{lemma}\label{BDGfaible}
Let $(M^n_{t \in [0,T]}) $ be a sequence of $(\shf_{t})$ local martingales such that
$M^n_0 = 0$ and
$<M^n>_T$ converges to $0$ in probability as $n \rightarrow \infty$.
Then $M^n  \rightarrow 0$ ucp, when $n \rightarrow \infty$.

\end{lemma}
\begin{proof}

It suffices to apply to $N=M^n$
the following inequality stated in \cite{ks}, Problem 5.25 Chap. 1,
which holds for any $(\shf_t)$-local martingale $(N_t)$ such that $N_0 = 0$:
\begin{equation}\label{SIVR32}
    P\big(\sup_{0\leq u \leq t}|N_u|\geq \lambda \big)\leq
   P\big(<N>_t \geq \delta \big)+
   \frac{1}{\lambda^2} E\big[\delta \wedge <N>_t\big],
\end{equation}
for any $t\geq 0, \ \lambda, \delta >0$.
\end{proof}
An immediate consequence of the  previous lemma is the following.

\begin{corollary} \label{C515} Let $M$ be an $(\shf_t)$-local martingale
vanishing at zero, with ${<M>} = 0$. Then $ M$ is identically zero.
\end{corollary}


\subsection{The It\^o integral} \label{SITOInt}



Let $M$ be an $(\shf_t)$-local martingale.
We construct here the It\^o integral with respect to $M$
using stochastic calculus via regularization.
We will proceed in two steps. First we  define the It\^o integral
$\int_0^\cdot  H dM$ for  a smooth integrand process $H$ 
as the forward integral $\int_0^\cdot  H d^-M$.
 Second, we extend $H \mapsto \int_0^\cdot  H dM$
via functional analytical arguments. We remark that the
classical theory of It\^o integrals first defines the integral of
simple  step processes $H$, see Remark \ref{RITO}, for details.

Observe first that the forward integral of
a continuous process $H$ of  bounded variation  is well defined because  Proposition \ref{P24}
4), 7)  imply that
\begin{equation} \label{SIVR32bis}
 \int_0^t Hd^-M = H_t M_t - H_0 M_0 - 
 \int_0^t M d^+H =  H_t M_t - H_0 M_0 - \int_0^t M_s dH_s.
\end{equation}
Call  $\shc$  the vector algebra of
 adapted processes
whose paths are of class  $C^0$.
This linear space,
equipped with the metrizable topology which governs the ucp convergence, is an $F$-space. For the definition
and properties of $F$-spaces, see  \cite{ds}, chapter 2.1.
Remark that the set  $\shm_{\rm loc}$ of 
continuous $(\shf_t)$-local martingales is
a closed  linear subspace of $\shc$, see for instance  \cite{gr}.

Denote by $\shc^{BV}$  the $\shc$ subspace of processes
whose paths are a.s. continuous with bounded variation.
The next observation is crucial.

\begin{lemma} \label{Contito}
If $H$ is an adapted process  in $\shc^{BV}$ then $\left (
\int_0^\cdot Hd^-M \right )$ is an  $(\shf_t)$-local martingale
whose quadratic variation is given by
$$ < \int_0^\cdot Hd^-M >_t =   < \int_0^\cdot H^2_s d <M>_s.$$
\end{lemma}
\begin{proof}

We only sketch  the proof. We restrict ourselves to  prove that if
$M$ is a local martingale then $\dis  Y=\int_0^\cdot Hd^-M $ is a
local martingale.

By localization, we can suppose that $H$, its total variation 
$\Vert H \Vert$
 and $M$ are bounded
processes.

Let $0 \leq s <t$. Since $H_t = H_0 + \int_s^t dH_u$,
  (\ref{SIVR32bis}) implies
\begin{equation}\label{SIVR32h}
    Y_t=H_sM_t -H_0M_0-\int_0^s M_u dH_u +\int_s ^t
(M_t-M_u)  dH_u.
\end{equation}
Let $(\pi_n)$ be a sequence of subdivisions of $[s,t]$, such that the mesh
of $(\pi_n)$ goes to zero when $n \rightarrow  + \infty$.
Since $M$ is continuous, $M$ and $\Vert H \Vert$ are bounded,
$$ \Delta_n := \sum_{\pi_n} (M_t - M_{u_{i+1}}) (H_{u_{i+1}} - H_{u_i}), $$
goes to $\int_s^t (M_t - M_u) dH_u$ a.s. and in $L^1$.
Consequently,
$$ E\left(\int_s^t (M_t - M_u) dH_u \right) = \lim_{n \rightarrow \infty} E(\Delta_n \vert \shf_s)$$
and
$$ E(\Delta_n \vert \shf_s) = \sum_{\pi_n} E \left((M_t - M_{u_{i+1}}) (H_{u_{i+1}} - H_{u_i}) \vert \shf_s \right ). $$
But one has
\begin{eqnarray}
  E \left((M_t - M_{u_{i+1}}) (H_{u_{i+1}} - H_{u_i}) 
  \vert \shf_s \right )\hskip-33mm &&\nonumber\\
  &=&  E \left( E( (M_t - M_{u_{i+1}}) (H_{u_{i+1}} - H_{u_i}) \vert \shf_{u_{i+1}} )
 \vert \shf_s \right ) \\
& = &  E \left((H_{u_{i+1}} - H_{u_i})  E(M_t - M_{u_{i+1}} 
\vert \shf_{u_{i+1}} ) \vert  \shf_s  \right)  \nonumber\\
& = & 0,
\end{eqnarray}
since $H$ is adapted and $M$ is a martingale.

Finally, taking the conditional expectation with respect to
${\cal F}_s$ in (\ref{SIVR32h}) yields
$$E\big[Y_t|{\cal F}_s\big]=H_sM_s -H_0M_0 -\int_0^s
M_u dH_u = Y_s.$$
Similar arguments show that $\dis Y^2-\int_0^\cdot
H^2d<M>$ is a martingale.

\end{proof}

The previous lemma  allows to extend the map $H \mapsto \int_0^t Hd^-M$.
Let   ${\shl^2(d<M>)}$ denote the set of  progressively measurable processes
such that
\begin{equation} \label{SIVR32prog}
 \int_0^T H^2 d<M> < \infty \ {\rm  a.s.}
\end{equation}
$\shl^2(d<M>)$ is an $F$-space with respect to the metrizable topology $d_2$
defined as follows:
  $(H^n)$ converges to $H$ when $n \rightarrow \infty$ if
$\int_0^T (H^n_s - H_s)^2 d<M>_s \rightarrow 0$ in probability, when $n \rightarrow \infty$.

\begin{remark} \label{RContIto}
$\shc^{BV}$ is dense in  $\shl^2(d<M>)$.
Indeed, according to \cite{ks}, lemma 2.7 section 3.2, simple processes are dense into  $\shl^2(d<M>)$.
On the other hand, a simple process of the form $H_t = \xi 1_{]a,b]}$, $\xi$ being $\shf_a$
measurable, can be expressed as a limit of  $H^n_t = \xi \phi^n$ where $\phi^n$ are  continuous
functions with bounded variation.

\end{remark}

Let $ \Lambda : \shc^{BV} \rightarrow \shm_{\rm loc}$
 be the map defined by
$\Lambda  H = \int_0^\cdot H d^-M$.

\begin{lemma}\label{Itoext}
If $\shc^{BV}$  (resp.  $\shm_{\rm loc}$) is equipped with $d_2$ (resp. the ucp topology)
then $ \Lambda$ is continuous.
\end{lemma}
\begin{proof}

Let $H^k$ be a sequence of processes in $\shc^{BV}$, converging 
to $0$ for $d_2$ when $k \rightarrow \infty$.
Set $N^k = \int_0^\cdot H^k d^-M$.
Lemma \ref{Contito} implies that $<N^k>_T$ converges to~$0$ in probability.
Finally Lemma \ref{BDGfaible} concludes the proof.
\end{proof}
We can now easily define the It\^o integral.
Since $\shc^{BV}$ is dense in $\shl^2(d<M>)$ for $d_2$,
Lemma \ref{Itoext} and standard functional analysis arguments imply that $\Lambda$ uniquely and continuously extends to  $\shl^2(d<M>)$.
\begin{definition} \label{Dito}
If $H$ belongs to  $\shl^2(d<M>)$,
we put $\int_0^\cdot H dM := \Lambda H$ and we call this the
{\bf It\^o integral of $H$ with respect to $M$}.
\end{definition}

\begin{proposition} \label{PITOExt}
If $H$ belongs to  $\shl^2(d<M>)$,  then $(\int_0^\cdot H dM)$ is an
$(\shf_t)$-local martingale with bracket
\begin{equation} \label{SIVR32ter}
<\int_0^\cdot H dM> = \int_0^\cdot H^2 d<M>.
\end{equation}
\end{proposition}

\begin{proof} \ Let $H\in \shl^2(d<M>)$. From Definition
\ref{Dito}, $(\int_0^\cdot H dM)$ is an $(\shf_t)$-local
martingale. It remains to prove  (\ref{SIVR32ter}).

Since $H$ belongs to $\shl^2(d<M>)$, then there exists a sequence
$(H_n)$ of elements in $\shc^{BV}$, such that $H_n\rightarrow H$ in
$\shl^2(d<M>)$.

Introduce $\dis N_n = \int_0^\cdot H_n dM$ and $\dis
N_n'=N_n^2-{<}N_n{>}$.
According to lemma~\ref{Itoext},
$\dis <N_n>=\int_0^\cdot H_n^2 d<M>$;
now  $N_n\rightarrow N$, ucp, $n\rightarrow \infty$ and
$ <N_n>  $
 goes to $\dis \int_0^\cdot H^2 d<M>$ in the ucp sense,
 as $n \rightarrow \infty$.
Therefore  $N'_n$ converges with respect to the ucp topology, to
the local martingale $\dis N^2-\int_0^\cdot H^2d<M>$. 
This actually proves (\ref{SIVR32ter}).

\end{proof}

\begin{remark}\label{RITO}
\begin{enumerate}
\item Recall that whenever  $H \in \shc^{BV} $
$$ \int_0^\cdot H dM = \int_0^\cdot H d^-M.$$
This property will be generalized 
in Propositions \ref{PITOV3bis} and \ref{TITOV2}.
\item
We emphasize  that It\^o stochastic 
integration based on adapted simple step processes
and the previous construction, finally lead to the same object.

If $H$ is of  the type
$ Y 1_{]a,b]} $ where $Y$ is an $\shf_a$ measurable random variable,
it is easy to show that $\int_0^t H dM = Y (M_{t \wedge b} - M_{t \wedge a})$.
Since the class of elementary processes obtained by linear combination
of previous processes is dense in $\shl^2(d<M>)$ and the map $\Lambda$ is continuous, then    $\int_0^\cdot H dM$
equals the classical It\^o integral.
\end{enumerate}
\end{remark}

In  Proposition
\ref{PITOV4} below we state the chain rule property.

\begin{proposition}\label{PITOV4}
 Let $(M_t, t \geq 0)$ be an  $({\cal F}_t)$-local
martingale, $(H_t, t \geq 0)$ be
in  $\shl^2(d <M>)$,
 $\displaystyle N :=
\int _0^\cdot H_sdM_s$ and $(K_t, t \geq 0)$ be a $({\cal
F}_t)$-progressively measurable process such that  $
\displaystyle \int _0^T (H_sK_s)^2d <M>_s<\infty$ a.s. Then
\begin{equation}\label{ITOV22}
\int_0^tK_sdN_s = \int _0^tH_sK_sdM_s, \quad 0 \leq t \leq T.
\end{equation}

\end{proposition}

\begin{proof} \ Since the map $\Lambda :  H\in {\cal
L}^2(d<M>)\mapsto \int_0^\cdot HdM$ is continuous, it
suffices to prove (\ref{ITOV22}) for $H$ and $K$ 
continuous and with bounded variation.

For simplicity we suppose $M_0=H_0=K_0=0$.

One has
$$\int_0^t KdN=\int_0^t(N_t-N_u)dK_u,$$
and
$$\begin{array}{ccl}
  N_t-N_u & = & \dis \int_0^t(M_t-M_v) dH_v-\int_0^u(M_u-M_v) dH_v \\
   & = & \dis (M_t-M_u)H_u +\int_u^t(M_t-M_v) dH_v,\\
\end{array}
$$
where $0 \le u \le t$.

Using Fubini's theorem one gets
\begin{eqnarray}
 \int_0^t KdN = \int_0^t (M_t-M_u)(H_udK_u + K_u dH_u)
              \hskip-33mm &&\nonumber\\
  &=&\int_0^t (M_t-M_u)  d (HK)_u=\int_0^t HKdM.\nonumber
\end{eqnarray}

\end{proof}

\subsection{Connections with calculus via 
regularizations} \label{SSConnReg}

The next Proposition  will show that, under suitable
conditions, the  It\^o integral is a
forward integral.

\begin{proposition}\label{PITOV3bis}
 Let $X$ be an $({\cal F}_t)$-local martingale
 and suppose that $(H_t)$  is progressively measurable and locally bounded.
\begin{enumerate}
\item If $H$
has a left limit at each point
then
 $\displaystyle \int_0^\cdot H_sd^-X_s=\int_0^\cdot
H_{s-} dX_s$.
\item If  $H_t=H_{t-}, \ d<X>_t$ a.e. (in particular if $H$ is c\`adl\`ag),
then $\displaystyle \int_0^\cdot H_sd^-X_s=\int_0^\cdot
H_{s} dX_s$.
\end{enumerate}
\end{proposition}

\begin{proof} \

Since $\dis s\mapsto \int _{s-\varepsilon}^s H_udu$ is continuous with bounded
variation,
$$\begin{array}{ccl}\dis
  \int_0^t\Big(\frac{1}{\varepsilon}\int _{s-\varepsilon}^s
H_udu\Big)d X_s & = & \dis\int_0^t\Big(\frac{1}{\varepsilon}\int
_{s-\varepsilon}^s
H_udu\Big)d^-X_s \\
&&\\
   & = &  \dis X_t\Big(\frac{1}{\varepsilon}\int
_{t-\varepsilon}^t H_udu\Big) - H_0 X_0  -\frac{1}{\varepsilon}
\int_0^t(H_s-H_{s-\varepsilon})X_sds.\\
\end{array}
$$
The second integral in the right-hand side can be modified
as follows
\begin{eqnarray*}
-\int_0^t(H_s-H_{s-\varepsilon})X_sds &= &\int_0^tH_s(X_{s+\varepsilon}-X_s)ds
-\int_{t-\varepsilon}^tH_sX_{s+\varepsilon}ds \\
&+& H_0 \int_0^\varepsilon X_s ds.
\end{eqnarray*}
Consequently
\begin{equation}\label{ITOV16}
    \int_0^t\Big(\frac{1}{\varepsilon}\int _{s-\varepsilon}^s
H_udu\Big)dX_s=\frac{1}{\varepsilon}\int_0^tH_s(X_{s+\varepsilon}-X_s)ds+R_\varepsilon(t),
\end{equation}
where
\begin{eqnarray}\label{ITOV15}
R_\varepsilon(t) & = & X_t\Big(\frac{1}{\varepsilon}\int
_{t-\varepsilon}^t
H_s ds\Big)-\frac{1}{\varepsilon}\int_{t-\varepsilon}^tH_sX_{s+\varepsilon}ds +  H_0 \left (\frac{1}{\varepsilon}
\int_0^\varepsilon X_s ds - X_0 \right)  \nonumber \\
&& \\
& = & \frac{1}{\varepsilon}\int_{t-\varepsilon}^t H_s (X_t -  X_{s+\varepsilon}) ds +  H_0 \left (\frac{1}{\varepsilon}
\int_0^\varepsilon X_s ds - X_0 \right) \nonumber
\end{eqnarray}
converges to zero ucp.

Under assumption 1,  Lebesgue's  dominated 
convergence theorem implies that
 $\frac{1}{\varepsilon}
 \int_{\cdot-\varepsilon}^\cdot H_s ds$  converges to $H_-$
according to ${\cal L}^2(d<M>)$, so the 
left-hand side of equality  (\ref{ITOV16})
converges to the  It\^o integral 
$\dis \int_0^\cdot H_{s-} dX_s$. This forces
the right-hand side  to converge to  
$\displaystyle \int_0^\cdot H_sd^-X_s$.

The proof of 2 is similar, remarking that
$H_s = H_{s-}$, for $d<M>_s$ a.e.\end{proof}

When the integrator is a Brownian motion $W$, we will see
in Theorem \ref{TITOV2} below
that the forward integral coincides  with the It\^o integral
for any integrand  in $\shl^2(d<W>).$
This is no longer true  when the integrator is
a general semimartingale.
The following example provides
a martingale $(M_t)$ and
a deterministic integrand $h$ such that 
the It\^o integral $\displaystyle
\int_0^t h dM$ and the forward integral $\dis \int_0^t h d^-M$ exist,
but are different.

\begin{example}\label{EITOV7} Let  $\psi :
[0,\infty[\mapsto\mathbb{R}$ verify $\psi (0)=0$, $\psi$ is
continuous, increasing, and $\psi'(t)=0$ a.e. (with respect to the
Lebesgue measure). Let $(M_t)$ be the process: $M_t= W_{\psi(t)},\
t\geq 0$, and $h$  be the indicator function of the support of the
positive measure $d\psi$.
Since $W^2_t - t$ is a martingale, $<W>_t = t$.
Clearly $(M_t)$ is a  martingale
and
    $<M>_t=\psi (t), \ t\geq 0$.
Observe that
 $h = 0$ a.e. with respect to Lebesgue measure.
Then
$ \dis \int _0^\cdot h(s)
\frac{M(s+\varepsilon)-M(s)}{\varepsilon} ds = 0 $
and so $\dis \int_0^\cdot h d^-M = 0.$

On the other hand,  $h=1, \ d\psi$ a.e., implies $\dis
\int_0^t h dM = M_t, \ t\geq 0$.
\end{example}

\begin{remark}\label{PITOV6}
A significant result of classical stochastic calculus is
the Bichteler-Dellacherie theorem, see \cite{prot} Th. 22, Section III.7.
In  the regularization approach,   an analogous property  occurs:
if the forward integral exists for a rich class
of adapted integrands, then the integrator is forced to be a
semimartingale. More precisely we recall the significant statement
of \cite{rv1}, Proposition 1.2.

 Let $(X_t, t\geq 0)$ be an
$({\cal F}_t)$-adapted and continuous process such that for any
c\`{a}dl\`{a}g, bounded and adapted
 process $(H_t)$, the forward integral $\displaystyle \int _0^\cdot
Hd^-X$ exists. Then $(X_t)$ is an $({\cal F}_t)$-semimartingale.
\end{remark}

From Proposition \ref{PITOV3bis} we deduce the relation between skew
and square bracket.

\begin{corollary}\label{CDrob}
Let $M$ be an  $(\shf_t)$-local martingale. Then
$<M> = [M]$ and
\begin{equation} \label{ESqobl}
M_t^2 = M_0^2 + 2 \int_0^t M d^-M + <M>_t .
\end{equation}
\end{corollary}

%
%
%

\begin{proof}\ The proof of (\ref{ESqobl}) is very simple and is based
on the following identity
$$(M_{s+\varepsilon}-M_s)^2=M_{s+\varepsilon}^2-M_s^2-2M_s(M_{s+\varepsilon}-M_s).$$
Integrating on $[0,t]$ leads to
$$\begin{array}{ccl}
 \dis  \frac{1}{\varepsilon}\int_0^t (M_{s+\varepsilon}-M_s)^2 ds& = &
  \dis  \frac{1}{\varepsilon}\int_0^t M_{s+\varepsilon}^2 ds
  -\frac{1}{\varepsilon}\int_0^t M_{s}^2 ds
-\frac{2}{\varepsilon}\int_0^tM_s(M_{s+\varepsilon}-M_s)ds
  \\
   & = &  \dis  \frac{1}{\varepsilon}\int_t^{t+\varepsilon} M_{s}^2 ds
  -\frac{1}{\varepsilon}\int_0^\varepsilon M_{s}^2 ds
-\frac{2}{\varepsilon}\int_0^tM_s(M_{s+\varepsilon}-M_s)ds.\\
\end{array}
$$
Therefore, taking the limit when $\varepsilon\rightarrow 0$, one
obtains
$$[M]_t=M_t^2-M_0^2-2\int_0^tM_sd^-M_s.$$
Since $t\mapsto M_t$ is continuous, the forward integral $\dis
\int_0^\cdot M d^-M$ coincides with the corresponding It\^{o}
integral. Consequently $M_t^2-M_0^2-[M]_t$ is a local martingale.
This proves both $[M]=<M>$ and (\ref{ESqobl}).
\end{proof}

\begin{corollary}\label{C1V}
Let $M,M'$ be two   $(\shf_t)$-local martingales.
Then $(M,M')$ has all its mutual covariations.
\end{corollary}

\begin{proof}
\ Since $M, M'$ and $M + M'$ are continuous local
martingales, Corollary~\ref{CDrob}
 directly implies that they
 have finite quadratic variation. The
bilinearity property of the covariation directly implies that
$[M,M']$ exists  and equals $$ \frac{1}{2} ([M + M'] - [M] -
[M']).$$
\end{proof}

\begin{proposition}\label{Crintst} Let $M$ and $M'$ be two
$(\shf_t)$-local martingales, $H$ and $H'$ be two progressively
measurable processes such that
$$\int_0^\cdot H^2d<M><\infty, \quad \int_0^\cdot
H^2d<M'><\infty.$$
Then
$${\Bigl[\int_0^\cdot H dM, \int_0^\cdot H' dM'\Bigr]}_t
=\int_0^tHH'd [M,M']_t.$$
\end{proposition}

The next proposition provides a simple example of two processes
$(M_t)$ and $(Y_t)$ such that $[M,Y]$ exists even though the
vector $(M,Y)$ has no mutual covariation.

\begin{proposition} \label{P611bis}
Let  $(M_t)$ be an continuous $({\cal F}_t)$-local martingale,
 $(Y_t)$  a c\`{a}dl\`{a}g  and an $({\cal F}_t)$-adapted  process.
If $M$ and $Y$ are independent then $[M,Y] = 0$.
\end{proposition}

\begin{proof}
\ Let $\shy$  be the $\sigma$-field generated by $(Y_t)$,
and denote by $(\tilde \shm_t)$ the smallest filtration satisfying 
the usual conditions
 and containing $( \shf_t)$ and $\shy$, i.e.,
 $\sigma(M_s, s \le t) \vee \shy \subset
 \tilde \shm_t, \forall t \geq 0$.
It is not difficult to show that $(M_t)$ is also an $(\tilde
\shm_t)$-martingale.

Thanks to Proposition   \ref{P24} 1., it is sufficient  to prove that
\begin{equation}\label{611}
    \int_0^t Y d^- M = \int_0^t Y d^+M.
\end{equation}
 Proposition \ref{PITOV3bis}  implies that the  left-hand side
 coincides with the $(\shm_t)$-It\^o integral $\int_0^t Y dM $.

\ Without restricting generality  we suppose
 $M_0=0$. We proceed as in the proof of Proposition
\ref{PITOV3bis}. Since a.s.  $\dis s\mapsto \int^ {s+\varepsilon}_s
Y_udu$ is continuous with bounded variation, 
$$
  \int_0^t\Big(\frac{1}{\varepsilon}\int ^ {s+\varepsilon}_s
Y_udu\Big)d ^-M_s = M_t\Big(\frac{1}{\varepsilon}\int ^
{s+\varepsilon}_sY_udu\Big)-\frac{1}{\varepsilon}
\int_0^t(Y_{s+\varepsilon}-Y_s)M_sds.
$$
As the processes $Y$ and $M$ are independent, the
forward integral in the left-hand side above is actually an
It\^{o} integral. Therefore, taking the limit 
when $\varepsilon
\rightarrow 0$ and using Proposition \ref{PITOV3bis}, one gets
$$ \int_0^tYdM=\int_0^tYd^-M=Y_tM_t-\int_0^tMd^-Y.$$
According to point 4) of Proposition \ref{P24}, the right-hand side
is equal to $\dis \int_0^tYd^+M$; this proves
(\ref{611}).
\end{proof}

\subsection{The semimartingale case} \label{SS54}

We begin this section with a technical lemma which implies
that the decomposition of a semimartingale is unique.

\begin{lemma} \label{lss54} Let $(M_t, t\geq 0 )$ be a
 $({\cal F}_t)$-local martingale with bounded variation. Then
 $(M_t)$ is constant.
\end{lemma}

\begin{proof} \ Since $M$ has bounded variation, then Proposition
\ref{P24}, 7) implies that $[M]=0$. Consequently Corollaries
\ref{C515} and \ref{CDrob}  imply that $M_t=M_0, \ t\geq 0$.

\end{proof}

It is now easy to define stochastic integration with respect to
continuous semimartingales.

\begin{definition} \label{DITOV2} Let $(X_t, t\geq 0 )$ be an
 $({\cal F}_t)$-semimartingale  with canonical
decomposition $X =M+V$, where $M$ (resp. $V$) is a continuous
$({\cal F}_t)$-local  martingale (resp. bounded variation,
continuous and $({\cal F}_t)$-adapted process) vanishing at $0$.
Let $(H_t, t\geq 0 )$ be  an $({\cal F}_t)$-progressively
measurable process, satisfying
\begin{equation}\label{ITOV17a}
\int_0^T H_s^2d[M,M]_s <\infty , \quad \mbox{and } \quad \int_0^T
|H_s| d\Vert V \Vert_s <\infty,
\end{equation}
where $\Vert V \Vert_t$ is the total variation of $V$ over $[0,t]$.

We set
$$
\int_0^t H_sdX_s=\int_0^tH_sdM_s+\int_0^tH_sdV_s, \quad 0\leq t
\leq T.$$
\end{definition}

\begin{remark} \label{RITOV3} \begin{enumerate}
\item In the previous definition, 
the integral with respect to $M$ (resp.
$V$) is an It\^{o}-type (resp. Stieltjes-type) integral.

\item It is clear that $\displaystyle \int_0^\cdot H_sdX_s$ is
again a continuous $({\cal F}_t)$-semimartingale, with martingale
part $\displaystyle \int_0^\cdot H_sdM_s$ and bounded variation
component $\displaystyle \int_0^\cdot H_sdV_s$.
\end{enumerate}
\end{remark}

Once we have introduced stochastic integrals with respect to
continuous semimartingales, it is easy to define Stratonovich
integrals.

\begin{definition} \label{DITOV3} Let $(X_t, t\geq 0 )$ be an
 $({\cal F}_t)$-semimartingale and $(Y_t, t\geq 0 )$
an $({\cal F}_t)$-progressively measurable process. The {\bf
Stratonovich} integral of $Y$ with respect to $X$ is defined as
follows
\begin{equation}\label{ITOV18}
\int_0^t Y_s\circ dX_s=\int_0^t Y_s dX_s +\frac{1}{2}[Y,X]_t;
\quad  t \geq 0,
\end{equation}
if $[Y,X]$ and $\int_0^\cdot Y_s dX_s$ exist.
\end{definition}

\begin{remark}\label{RITOV4} \begin{enumerate}
\item Recall that conditions of  type (\ref{ITOV17a}) ensure
existence of the stochastic integral with respect to  $X$.

\item If $(X_t) $ and $(Y_t)$ are   $({\cal
F}_t)$-semimartingales, then $\int_0^\cdot Y_s\circ dX_s$ exists
and is called the {\bf Fisk-Stratonovich} integral.

\item Suppose that $(X_t  )$ is an $({\cal
F}_t)$-semimartingale and $(Y_t  )$ is a left continuous and $({\cal
F}_t)$-adapted  process such that $[Y,X]$ exists. We already have
observed (see Proposition \ref{PITOV3bis}) that $\displaystyle
\int_0^\cdot Y_s dX_s$ coincides with $\displaystyle \int_0^\cdot
Y_s d^-X_s$. Proposition \ref{P24} 1) and 2) imply that the
Stratonovich integral $\displaystyle \int_0^\cdot Y_s\circ dX_s$
is equal to the symmetric integral $\displaystyle\int_0^\cdot Y_s
d^\circ X_s$.

\end{enumerate}
\end{remark}

At this point we can easily identify the  covariation of two
semimartingales.

\begin{proposition} \label{QVS}
Let $S^i = M^i + V^i$ be two $(\shf_t)$-semimartingales, $i=1,2$,
where $M^i$ are local martingales and $V^i$ bounded variation
processes. One has $[S^1,S^2] = [M^1,M^2]$.
\end{proposition}

\begin{proof} \ The result follows directly from
Corollary \ref{C1V}, Proposition \ref{P24}  7), and the
bilinearity of the covariation.
\end{proof}
\begin{corollary}\label{C617}
Let $S^1, S^2$ be two  $(\shf_t)$- semimartingales such that their
martingale parts are independent. Then $[S^1, S^2] = 0.$
\end{corollary}
\begin{proof}

It follows from Proposition \ref{P611bis}.
\end{proof}

The statement of
Proposition \ref{PITOV3bis} can be adapted to 
semimartingale integrators as follows.

\begin{proposition} \label{PITOV7} Let 
$X$ be an $({\cal F}_t)$-semimartingale
 and suppose that $(H_t)$  is adapted,
with left limits at each point.
 Then
 $\displaystyle \int_0^\cdot H_sd^-X_s=\int_0^\cdot
H_{s-} dX_s$.
If $H$ is c\`adl\`ag then
 $\displaystyle \int_0^\cdot Hd^-X=\int_0^\cdot
H  dX$.
\end{proposition}

\begin{remark} \label{Rpatper}
\begin{enumerate}
\item Forward integrals  generalize not only classical It\^o
integrals but also the integral obtained from the
theory of  enlargements of
filtrations, see e.g. \cite{jy}. Let $(\shf_t)$ and
$(\shg_t)$ be two filtrations fulfilling 
the usual conditions with
$\shf_t \subset \shg_t$ for all $t$. Let $X$ be a
$(\shg_t)$-semimartingale which is $(\shf_t)$-adapted.
By Stricker's theorem, $X$ is also an $(\shf_t)$-semimartingale.
Let $H$ be a c\`adl\`ag  bounded $(\shf_t)$-adapted process.
According to Proposition \ref{PITOV7}, the $(\shf_t)$-It\^o
integral $\int_0^\cdot HdX$ equals the  $(\shg_t)$-It\^o integral
and it coincides with the forward integral  $\int_0^\cdot Hd^- X$.
\item The result stated above is false when $H$ has no left limits
at each point. Using a tricky example in \cite{perkins}, it is
possible to exhibit a filtration $(\shg_t)$, a $(\shg_t)$-semimartingale ${(X_t)}_{t \ge 0}$
with natural filtration $\shf^X_t$,
 a bounded and $(\shf^X_t)$-progressively measurable process $H$, such that
 $\int_0^\cdot H d^-X$ equals the $(\shf^X_t)$-It\^o integral
but  differs from the $(\shg_t)$-It\^o integral. More precisely one
has:
\begin{enumerate}
\item $X$ is a  3-dimensional Bessel process
with decomposition
\begin{equation} \label{EBessel}
 X_t = W_t + \int_0^t \frac{1}{X_s} ds,
\end{equation}
where $W$ is an $(\shf^X_t)$-Brownian motion,
 \item  $X$ is a
$(\shg_t)$-semimartingale with  decomposition $M + V$
where $M$ is the local martingale part,
\item $H_t(\omega) = 1$ for $ dt {\otimes }dP$-almost all  
$(t,\omega) \in [0,T]
\times \Omega$, \item  $\beta_t = \int_0^t H dX $ is a
$(\shg_t)$-Brownian motion.
\end{enumerate}
Property (d) implies that $ I^-(\varepsilon,H, dX) =
I^-(\varepsilon,1 ,dX) $ so that $\int_0^t H d^- X = X_t$. The
$(\shf^X_t)$-It\^o integral $\int_0^t H dX $ equals $ \int_0^t H
dW +  \int_0^t \frac{H_s}{X_s} ds$; Theorem~\ref{TITOV2} below
and Proposition \ref{P24} 8) imply that this integral
coincides with $\int_0^t H d^- X$. Since a Bessel process cannot
be equal to a Brownian motion, the $(\shg_t)$-It\^o integral
$\int_0^t H d X$ differs from the $(\shf^X_t)$-It\^o integral
$\int_0^t H d X$.

Indeed, the pathology comes from the integration 
with respect to the bounded variation process.
In fact, according to ii), $[X]_t =  [W]_t = t$; therefore
   $M$ is a $(\shg_t)$-Brownian motion.
Theorem \ref{TITOV2} below says that $\int_0^\cdot H d^-M =
\int_0^\cdot H dM$;
 the additivity of forward integrals and It\^o integrals
imply that $\int_0^\cdot H d^-V \neq \int_0^\cdot H dV$.
Consequently it can be deduced from Proposition \ref{P24} 7) a)
that the discontinuities of $H$ are not a.s. countable.
It can even be shown 
that the discontinuities of $H$ are not
negligible with respect to $dV$.

\end{enumerate}
\end{remark}

\subsection{The Brownian case} \label{SS54B}

In this section we will investigate the link between forward and
It\^{o} integration  with respect to a Brownian motion. In this
section $(W_t)$ will denote a  $({\cal F}_t)$-Brownian motion.

The main result of this subsection is the following.

\begin{theorem}\label{TITOV2}
Let  $(H_t,  t \geq 0)$ be an $({\cal F}_t)$-progressively
measurable process satisfying $\displaystyle \int _0
^TH_s^2ds<\infty$ a.s. Then the It\^{o} integral $\displaystyle
\int_0^\cdot H_sdW_s$ coincides with the forward integral
$\displaystyle \int_0^\cdot H_sd^-W_s$.
\end{theorem}

\begin{remark}\label{RITOV2}
\begin{enumerate}
\item
 We would like to  illustrate the
advantage of using regularization instead of discretization
(\cite{fo}) through the following example.

Let $g$ be the indicator function of $\mathbb{Q}\cap
\mathbb{R}_+$.

Let $\Pi=\{t_0=0,t_1, \cdots ,t_N=T\}$ be a subdivision of $[0,T]$
and
$$I(\Pi,g,dW)_t:=\sum_ i g(t_i )\big(W(t_{i+1}\wedge t)-
W(t_{i}\wedge t)\Big); \quad 0 \leq t \leq T.$$

We remark that
$$I(\Pi,g,dW)_t=
\left\{
\begin{array}{cl}
  0 & \mbox{if} \ \Pi \subset \mathbb{R} \setminus \mathbb{Q}\\
  W_t & \mbox{if} \ \Pi \subset  \mathbb{Q}. \\
\end{array}
\right.
$$
Therefore there is no canonical definition of $\displaystyle
\int_0 ^t gdW$ through discretization. This is not surprising since
$g$ is not a.e. continuous and so is not Riemann integrable. On
the contrary, integration via regularization seems drastically more
adapted to define $\displaystyle \int_0 ^t gd^-W$, for any $g \in
L^2([0,T])$, since this integral coincides with the classical
It\^{o}-Wiener integral.
\item
In order to overcome this problem, McShane pointed  
out an alternative approximation scheme, 
  see \cite{mcshane} chap. 2 and 3.
McShane's stochastic integration  makes use of the 
so-called {\it belated} partition;
the integral is then even more general than It\^o's one,
and
it includes in particular
the function $g$ above.

\end{enumerate}
\end{remark}

\begin{proof} \ (of  Theorem \ref{TITOV2}) \ 1) First,
suppose  in addition that $H$ is
 a continuous process. Replacing
$X$ by $W$ in (\ref{ITOV16}) one gets

\begin{equation}\label{ITOV14}
    \int_0^t\Big(\frac{1}{\varepsilon}\int _{s-\varepsilon}^s
H_udu\Big)dW_s=\frac{1}{\varepsilon}\int_0^tH_s(W_{s+\varepsilon}-W_s)ds+R_\varepsilon(t),
\end{equation}
where the remainder term $R_\varepsilon(t)$ is given by
(\ref{ITOV15}).

Recall the   maximal  inequality (\cite{st70}, chap. I.1): there
exists
 a constant $C$ such that for any $\phi \in L^2([0,T])$,
\begin{equation}\label{stein}
\int_0^T\Big(\sup_{0<\eta <1}
\Big\{\frac{1}{\eta}\int^v_{(v-\eta)_+}
 \phi_vdv \Big\}\Big)^2du \leq C \int_0^T \phi_v^2dv .
\end{equation}
2) We claim that (\ref{ITOV14}) may be extended to
any progressively
measurable process $(H_t)$  satisfying $\displaystyle \int _0
^\cdot H_s^2ds<\infty$.

Set $\displaystyle H^n_t=n \int_{t-1/n}^tH_udu$ for $t\geq 0$.
It is clear that as
    $n\rightarrow\infty$
\begin{itemize}
    \item for a.e. $t$, $H^n_t$ converges to $H_t$,
    \item $(H^n_t)$ converges to $(H_t)$ in ${\cal L}^2(d<W>)$
    (i.e. $\dis \int_0^\cdot (H_s^n-H_s)^2ds$ goes to $0$ in the
    ucp sense).
\end{itemize}
 Since
 $$ < \int_0^\cdot\Big(\frac{1}{\varepsilon}\int
_{s-\varepsilon}^s
H_udu\Big)dW_s>_t=\int_0^\cdot\Big(\frac{1}{\varepsilon}\int
_{s-\varepsilon}^s H_udu\Big)^2ds,$$

(\ref{stein}) and Lemma \ref{BDGfaible} imply that
(\ref{ITOV14}) and (\ref{ITOV15}) are  valid.

3)  Letting $\varepsilon \rightarrow 0$ in (\ref{ITOV14})
and using once more (\ref{stein}), Lemma     \ref{BDGfaible}
 allows to conclude the proof of  Theorem \ref{TITOV2}.

\end{proof}

\subsection{Substitution formulae} \label{SSubsti}

We conclude  Section~\ref{S5} by observing that discretization 
makes it possible 
to integrate  non adapted integrands in a context which is
covered neither by Skorohod integration theory nor by enlargement
of filtrations. A  class of  examples is the following.

Let $(X(t,x), t \geq 0, x \in \R^d)$ and 
 $(Y(t,x)$, $t \geq 0$, $x \in \R^d)$
be  two  families of continuous $({\cal F}_t)$ semimartingales
depending on a parameter $x$ and $(H(t,x)$, $t \geq 0$, $x \in \R^d)$
an $({\cal F}_t)$ progressively measurable processes depending on $x.$
 Let $Z$ be a ${\cal F}_T$-measurable r.v., taking its values in
$\mathbb{R}^d$.

Under some minimal  conditions of Garsia-Rodemich-Rumsey type,
 see for instance \cite{rv2, rv4}, one has
$$
\int_0^t H(s,Z)\, d^- X(s,Z) = \int_0^t H(s,x)\, dX (s,x)\Big|_{x =
Z},
$$
$$
[X(\cdot,Z), Y(\cdot,Z)] = [X(\cdot, x), Y(\cdot, x)] \Big|_{x = Z}.
$$
The first result is   useful to prove existence results  for SDEs
driven by semimartingales, with anticipating initial
conditions.

It is significant to remark that these  substitution formulae
give rise to anticipating calculus in a setting which is not covered
by Malliavin non-causal calculus since our integrators  may be
general semimartingales, while Skorohod integrals apply
essentially to Gaussian integrators or eventually to Poisson type
processes. Note that the usual causal It\^o calculus does
not apply here since $(X(s,Z))_s$ is not a semimartingale (take for
instance a r.v. $Z$ which generates  ${\cal F}_T$.)

\section{Calculus for finite quadratic variation processes}
\label{S6}

\subsection{Stability of the covariation} \label{S62}

 A basic tool of calculus via regularization 
is the stability of finite quadratic variation processes 
under
$C^1$ transformations.

\begin{proposition} \label{P32}
Let $(X^1, X^2)$ be a vector of processes having all its mutual
covariations and $f, g \in C^1(\R)$. Then $[f(X^1), g(X^2)]$ exists
and is given by
$$
[f(X^1),g(X^2)]_t = \int_0^t f'(X^1_s) g'(X^2_s) d[X^1,X^2]_s
$$
\end{proposition}

\begin{proof}
By polarization and bilinearity,
it suffices to consider the case
when $X = X^1 = X^2$ and $f = g$. 
Using Taylor's formula, one can write
$$
f(X_{s + \eps}) - f(X_s) = f'(X_s) (X_{s+\eps} - X_s) + R(s,\eps)
(X_{s+\eps} - X_s), \quad s \ge 0, \ \eps > 0,
$$
where $R(s,\varepsilon)$ denotes
a process which  converges
in the ucp sense to 0 when $\varepsilon\rightarrow 0$.
Since $f'$ is unifomly continuous on compacts,
$$
\left ( f(X_{s +\eps}) - f(X_s) \right )^2 = f'(X_s)^2 (X_{s +
\eps} - X_s)^2 +
  R(s,\eps) (X_{s +\eps}-X_s)^2.
$$
Integrating from $0$ to $t$ yields
$$
\frac{1}{\eps} \int_0^t (f(X_{s +\eps}) - f(X_s))^2 ds =
I_1(t,\eps) + I_2(t,\eps)
$$
where 
\beqnar
I_1(t,\eps) &=& \int_0^t f'(X_s)^2 \frac{(X_{s+\eps}-X_s)^2}{\eps} ds, \\
I_2(t,\eps) &=& \frac{1}{\eps} \int_0^t R(s,\eps)
(X_{s+\eps}-X_s)^2 ds. 
\eeqnar 
Clearly one has
$$
\sup_{t \le T} |I_2(t,\eps)| \le \sup_{s \le T}|R(s,\eps )|
\frac{1}{\eps} \int_0^T ( X_{s+\eps} - X_s)^2 ds.
$$
Since $[X]$ exists, $I_2(\cdot, \eps) \stackrel{\rm
ucp}{\longrightarrow} 0$. The result will follow if we establish
\begin{equation} \label{F34}
\frac{1}{\eps} \displaystyle \int_0^\cdot Y_s d\mu_\eps (s) \stackrel{\rm
ucp}{\longrightarrow}
 \int_0^\cdot Y_s
d[X,X]_s
\end{equation}
where $\mu_\eps (t) = \int_0^t \frac{ds}{\eps} (X_{s+\eps}-X_s)^2$
and $Y$ is a continuous process.
It is not difficult to verify that a.s., $ \mu_\varepsilon (dt)$
converges to $d[X,Y]$, when $\varepsilon \rightarrow 0$;
this finally implies (\ref{F34}).
\end{proof}

\subsection{It\^o formulae for
 finite quadratic variation processes} \label{S63}

Even though all It\^o formulae that we will consider can be stated
in the multidimensional case, see for instance \cite{rv2}, we
will only deal here  with dimension~1. Let $X = {(X_t)}_{t \ge 0}$
be a continuous process.

\begin{proposition} \label{P34}
Suppose that $[X,X]$ exists and let $f \in C^2(\rz{R})$. Then
\begin{equation} \label{F36}
\int_0^\cdot f'(X) d^-X \quad {\rm and} \quad \displaystyle \int_0^\cdot f'(X)
d^+X \quad {\rm
  exist}.
\end{equation}
Moreover
\begin{description}
\item{a)} $f(X_t) = f(X_0) + \int_0^t f'(X) d^\mp X \pm \halb
\int_0^t f''(X_s)
  d[X,X]_s$,
\item{b)} $f(X_t) = f(X_0) +  \displaystyle  \int_0^t f'(X) d^\mp X \pm \halb
[f'(X),X]_t$, \item{c)} $f(X_t) = f(X_0) + \displaystyle  \int_0^t f'(X) d^\circ
X $.
\end{description}
\end{proposition}

\begin{proof} \quad c) follows from b) summing up $+$ and $-$.

b) follows from a),  since Proposition  \ref{P32}  implies
that
$$
[f'(X), X]_t = \int_0^t f''(X) d[X,X].
$$
The proof of a) and (\ref{F36}) is similar to that of
 Proposition \ref{P32}, but with a second-order Taylor expansion.
\end{proof}

The next lemma emphasizes that the existence of a 
quadratic variation is
closely connected with the existence of some 
related forward and
backward integrals.

\begin{lemma} \label{L35}
Let $X$ be a continuous process. Then $[X,X]$ exists
$\Longleftrightarrow$
 $\displaystyle \int_0^\cdot X d^- X$ exists $\Longleftrightarrow$ $\displaystyle \int_0^\cdot X d^+X$ exists.
\end{lemma}

\begin{proof}

Start with the identity
\begin{equation} \label{F39}
(X_{s+\eps} - X_s)^2 = X_{s+\eps}^2 - X_s^2 - 2
X_s(X_{s+\eps}-X_s)
\end{equation}
and observe that, when $\eps \rightarrow 0,$
$$\frac{1}{\eps} \int_0^t
( X_{s+\eps}^2 - X_s^2) ds \rightarrow X^2_t - X^2_0.$$
Integrating (\ref{F39}) from $0$ to $t$ and dividing by $\eps$
easily gives the equivalence between the first two assertions.

The equivalence between the first  and third ones is similar,
  replacing $\eps$
with $ - \eps$ in (\ref{F39}).
\end{proof}

Lemma \ref{L35} admits the following generalization.

\begin{corollary} \label{C36}
Let $X$ be a continuous process. The following properties are
equivalent
\begin{description}
\item{a)} $[X,X]$ exists; \item{b)} $  \displaystyle 
\forall g \in C^1,\quad\int_0^\cdot g(X) d^-X$
exists; \item{c)} $ \displaystyle  
\forall g \in C^1,\quad\int_0^\cdot g(X) d^+X$
exists.
\end{description}
\end{corollary}

\begin{proof} \quad
The It\^o formula stated in Proposition \ref{P34} 1) implies a)
$\Rightarrow$ b). b) $\Rightarrow$ a) follows setting $g(x) =x$
and using Lemma \ref{L35}.

b) $\Leftrightarrow$ c) because of Proposition \ref{P24} 1) which
states that
$$
\int_0^\cdot g(X) d^+X = \int_0^\cdot g(X) d^-X + [g(X),X],
$$
and Proposition \ref{P32} saying that $[g(X),X]$
exists.
\end{proof}

When  $X$ is a
semimartingale, the It\^o formula 
seen above becomes the following.

\begin{proposition} \label{P611d} Let $(S_t)_{t \ge 0}$ be a continuous  $(\shf_t)$-semimartingale and  $f$
a function in $ C^{2} (\R)$.
One has the following.
\begin{enumerate}
\item
$$
f(S_t) = f(S_0) +
  \int_0^t f'(S_u) d S_u +
\frac{1}{2} \int_0^t f''(S_u) d[S,S]_u.
$$
\item Let $(S_t^0) $ be another continuous  $(\shf_t)$-semimartingale. The
following integration by parts holds:
$$S_t S^0_t = S_0 S^0_0 + \int_0^t S_u dS_u^0 +
 \int_0^t S^0_u dS_u + [S,S^0]_t.$$
\end{enumerate}
\end{proposition}

\begin{proof}

 We recall that It\^o and forward integrals coincide, see
Proposition \ref{PITOV3bis}; therefore point
1 is a consequence of Proposition \ref{P34}.

2  stems from the
integration by parts formula in Proposition \ref{P24} 4).

\end{proof}

\subsection{L\'evy area} \label{SSLevy}

In Corollary \ref{C36},  we have seen that  
$\int_0^t g(X) d^-X$ exists
when $X$ is a one-dimensional finite quadratic variation 
process and $g \in C^1(\R)$.

If $X = (X^1,X^2)$ is two-dimensional 
 and  has all its mutual covariations, 
consider $g \in C^1(\R^2; \R^2)$.
We naturally define, if it exists,
$$ \displaystyle \int_0^t g(X) \cdot  d^- X = 
\lim_{\varepsilon \rightarrow 0^+} 
I^-(\varepsilon , g(X) \cdot
dX)(t),$$
where
\begin{equation}\label{SIVR1vect}
I^-(\varepsilon , g(X) \cdot dX)(t)=\int _0^t g(X) (s)
\cdot  \frac{X(s+\varepsilon)-X(s)}{\varepsilon}ds; 
\quad 0 \leq t \leq
T,
\end{equation}
and $\cdot$ denotes the scalar product in $\R^2$.

With a  2-dimensional It\^o formula of the 
same type as in Proposition~\ref{P34},
it is possible to show that
  $\int_0^t g(X) \cdot  d^-X$
exists  if $g = \nabla u$, where $u$ is a potential of class $C^2$.
 If $g$ is a general $C^1(\R^2)$ function, 
 one cannot expect in general that $\int_0^t g(X)\cdot d^-X$
exists.

T. Lyons'  rough paths approach, 
see for instance \cite{lyons, lq, lejay, 
gubi, cq} has considered in detail
the problem of the existence of integrals of 
the type $\int_0^t g(X) \cdot dX $. In this 
theory, the concept of L\'evy area
plays a significant role.
Translating this in the present context one would say 
that the essential assumption is that $X = (X^1,X^2)$ has a
L\'evy area type process.
This section will only make some basic 
observations on that topic from
the perspective of stochastic calculus via regularization.

Given two classical semimartingales $S^1, S^2$, the classical
notion of L\'evy area is defined by
 $$ L(S^1, S^2)_t =  \int_0^t S^1 dS^2 -  \int_0^t S^2 dS^1,$$
 where both integrals are of It\^o type.

\begin{definition} \label{DLevy}
 Given two continuous processes $X$ and $Y$, we put
$$ L(X,Y)_t = \lim_{\varepsilon \rightarrow 0^+}
\int_0^t \frac{X_s Y_{s+\varepsilon} -  
X_{s+\varepsilon} Y_s}{\varepsilon}
ds.
$$
where the limit is understood in the ucp sense.
$L(X,Y)$ is called the {\bf L\'evy area} of 
the processes $X$ and $Y$.
\end{definition}
\begin{remark} \label{RLevy1}
The following properties are easy to establish.
\begin{enumerate}
\item   $L(X,X) \equiv 0$.
\item The L\'evy area is an antisymmetric operation, i.e.
$$L(X,Y) = - L(Y,X).$$
\end{enumerate}
\end{remark}
Using the approximation of symmetric integral
we  can easily prove the following.
\begin{proposition} \label{SLevy}
$\int_0^\cdot X d^\circ Y$ exists if and only if $L(X,Y)$ exists. Moreover
$$ 2 \int_0^t X d^\circ Y = X_t Y_t - X_0 Y_0 + L(X,Y)_t$$
\end{proposition}

Recalling  the convention
that an equality among three objects implies that at least 
two among the three are
defined, we have the following.

\begin{proposition} \label{PLevy}
\begin{enumerate}
\item $\dis \quad L(X,Y)_t = \int_0^t X d^\circ Y - \int_0^t Y d^\circ X.$
\item  $ \dis \quad L(X,Y)_t = \int_0^t X d^- Y - \int_0^t Y d^-X. $
\end{enumerate}
\end{proposition}
\begin{proof}

 \begin{enumerate}
\item
 From  Proposition \ref{SLevy} applied to $X,Y$ and $Y,X$,
and by antisymmetry of L\'evy areas we have
\begin{eqnarray*}
2 \int_0^t X d^\circ Y &=& X_t Y_t - X_0 Y_0 + L(X,Y)_t, \\
2 \int_0^t Y d^\circ X &=& X_t Y_t - X_0 Y_0 -  L(X,Y)_t.
\end{eqnarray*}
Taking the difference gives  1. 
\item  follows from the definition of forward integrals.
\end{enumerate}
\end{proof}

\begin{remark} \label{RLevy}
If $[X,Y]$ exists, point 2 of Proposition \ref{PLevy}  is a consequence
of point~1 and of  Proposition \ref{P24} 1, 2.
\end{remark}

For a real-valued process ${(X_t)}_{t \ge 0}$, Lemma
\ref{L35}  says that
$$ [X,X] \quad {\rm exists} \Leftrightarrow  \int_0^\cdot X d^-X \quad {\rm exists}.$$
 Given a vector of processes $\underline X =
(X^1, X^2)$  we may ask wether the following statement is true:

 $(X^1,X^2)$ has all its mutual
brackets if and only if
$$ \int_0^\cdot X^i d^-X^j \quad {\rm exists}, $$
for $i,j = 1,2$. In fact the answer is negative if the two-dimensional process $X$
does not have a L\'evy area.

\begin{remark} \label{RFQV} Suppose that $(X^1,X^2)$ has all its
mutual covariations. Let $ *$ stand for $\circ$, or $-$, or $+ $.
The following are equivalent.
\begin{enumerate}
\item The L\'evy area $L(X^1,X^2)$ exists.
\item $\int_0^\cdot X^i d^* X^j$ exists for any $i, j = 1,2$.
\end{enumerate}
By Lemma \ref{L35}, we first observe that $\int X^i d^\circ X^i$  exists
since $X^i $ is a finite quadratic variation process.
In  point 2, equivalence between the three cases
$ \circ$,$ -$ and $ + $ is obvious using  Proposition \ref{P24} 1 2.
Equivalence between the existence of $\int_0^\cdot X^1 d^\circ X^2$
and $L(X^1,X^2)$ was already established in Proposition~\ref{SLevy}.
\end{remark}

\section{Weak Dirichlet processes} \label{S8}

\subsection{Generalities} \label{S81}

Weak Dirichlet processes constitute a natural generalization
 of Dirichlet processes, which in turn
naturally extend semimartingales.
Dirichlet processes have been considered by many authors,
see for instance \cite{fodir, ber}.

Let $(\shf_t)_{t \ge 0}$ be a fixed  filtration
fulfilling the usual conditions.
In the present section~\ref{S8}, $(W_t)$ will 
denote a classical $(\shf_t)$-Brownian motion.
For simplicity, we shall stick to the framework 
of continuous processes.

\begin{definition} \label{D81}
\begin{enumerate}
\item An  $(\shf_t)$-{\bf Dirichlet process} is the sum of an
$(\shf_t)$-local martingale $M$ and a zero 
quadratic variation process $A$.
\item An  $(\shf_t)$-{\bf weak Dirichlet process}  
is the sum of an
$(\shf_t)$-local martingale $M$ and a process $A$ 
such that $[A,N] = 0$
for every continuous $(\shf_t)$- local martingale $N$.
\end{enumerate}
In both cases, we will suppose $A_0 = 0$ a.s.
\it
\end{definition}

\begin{remark}\label{R81}
\begin{enumerate}
\item  The process $(A_t)$ in the latter decomposition  is 
$(\shf_t)$-adapted.
\item Any $(\shf_t)$-semimartingale is an  
$(\shf_t)$-Dirichlet process.
\end{enumerate}
\end{remark}

The statement of the following  proposition is 
essentially contained in \cite{er2}.
\begin{proposition} \label{P81}
\begin{enumerate}
\item Any  $(\shf_t)$-Dirichlet process is an $(\shf_t)$-weak 
Dirichlet process.
\item The decomposition $M + A$ is unique.
\end{enumerate}
\end{proposition}
\begin{proof}

Point 1 follows from Proposition \ref{P24} 6).

Concerning point 2, let $X$ be a weak Dirichlet process with
decompositions
 $X = M^1 + A^1 = M^2 + A^2$. Then
$ 0 = M + A$  where $M = M^1 - M^2, A = A^1 - A^2$.
We evaluate the covariation of both members against $M$ to obtain
$$ 0 = [M] + [M,A^1] - [M,A^2] = [M].$$
Since $M_0 = A_0 = 0$ and $M$ is a local martingale, Corollary \ref{C515}
gives  $M = 0$.

\end{proof}

The class of semimartingales with respect to a given 
filtration is known
to be stable with respect to $C^2$ transformations, as
Proposition \ref{P611d} implies.
Proposition  \ref{P32} says  that finite quadratic variation processes
are stable under $C^1$ transformations.

It is possible to show that the class of weak Dirichlet processes
with finite quadratic variation
(as well as Dirichlet processes) is stable with respect to
the same type of transformations.
We start with a result which is a slight improvement 
(in the continuous case)
of a result obtained by \cite{cjms}.

\begin{proposition} \label{PDir}
Let $X$ be  a finite quadratic variation process 
which is $(\shf_t)$- weak Dirichlet, and
$f \in C^1(\R)$. Then $f(X)$ is also weak Dirichlet.
\end{proposition}
\begin{proof}

Let $X = M + A$ be the corresponding decomposition.
We express $f(X_t) =  M^f + A^f$
where
$$M^f_t = f(X_0)   +  \int_0^t f'(X) dM, 
\quad A^f_t = f(X_t) - M^f_t.$$
Let $N$ be a local martingale.
We have to show that $[f(X) - M^f, N] = 0.$

By additivity of the covariation,
and the definition of weak Dirichlet process,  $[X,N] = [M,N]$
so that Proposition \ref{P32} implies $[f(X), N]_t 
= \int_0^t f'(X_s) d[M,N]_s.$

On the other hand, Proposition \ref{Crintst}  gives
$$[ M^f, N]_t = \int_0^t   f'(X_s) d[M,N]_s, $$
and the result follows.
\end{proof}

\begin{remark} \label{RDir}
\begin{enumerate}
\item If $X$ is an $(\shf_t)$- Dirichlet process, 
it can be proved  similarly that
$f(X)$ is an $(\shf_t)$- Dirichlet process;
see  \cite{ber} and \cite{rvw} for details.
\item
The class of Lyons-Zheng processes introduced in \cite{rvw}
consitutes a natural generalization of reversible semimartingales, see Definition \ref{D89}.
The authors proved that  this class is also stable
through $C^1$ transformation.
\item Suppose that  $(\shf_t)$ is the canonical 
filtration associated with a Brownian motion $W$.
Then a continuous  $(\shf_t)$-adapted process $D$  
is weak Dirichlet if and only if $D$
is  it is the sum of an $(\shf_t)$-local martingale 
and a process $A$ such that $[A, W] = 0$.
See \cite{cr}, Corollary 3.10.
\end{enumerate}
\end{remark}

We  also report a Girsanov type theorem established by \cite{cjms}
at least in a discretization framework.

\begin{proposition} \label{Girsanov}
Let $X = {(X_t)}_{t \in [0,T]} $ be an  
$(\shf_t)$-weak Dirichlet process, and
 $Q$ a probability equivalent
 to $P$ on $\shf_T$.
Then $X = {(X_t)}_{t \in [0,T]} $
 is an $(\shf_t)$-weak Dirichlet process
with respect to $Q$.
\end{proposition}

\begin{proof}

We set $D_t = \frac{dQ}{dP} \vert_ {\shf_t}$;
$D$ is a positive local martingale.

Let $L$ be the local martingale such that 
$D_t = \exp(L_t - \frac{1}{2} [L]_t)$.
Let $X = M + A$ be the corresponding decomposition.
It is well-known that $\tilde M = M - [M,L]$ is a local martingale
under $Q$. So, $X$ is a $Q$-weak Dirichlet process.
\end{proof}

As mentioned earlier,  Dirichlet processes are stable with respect to
$C^1$ transformations. In applications, in particular to control theory,
one often needs to know the nature  of process $(u(t,D_t))$ where $u \in C^{0,1}(\R_+ \times \R)$
and $D$ is a Dirichlet process.
The following result was established in  \cite{gr}.

\begin{proposition} \label{Pgr}
Let $(S_t)$ be a continuous   $(\shf_t)$-weak Dirichlet process with finite quadratic variation; let  $u \in C^{0,1}(\R_+ \times \R)$.
Then $(u(t, S_t))$ is a  $(\shf_t)$-weak Dirichlet process.
\end{proposition}

\begin{remark} \label{Rgr}
There is no reason for  $(u(t,S_t))$ to have  
a finite quadratic variation 
since the dependence of $u$ on the first argument $t$ may be very rough.
A fortiori $(u(t,S_t))$ will not be Dirichlet.
Consider for instance $u$ only depending on time, deterministic,
with infinite quadratic variation.

\end{remark}

Examples of Dirichlet processes (respectively 
weak Dirichlet processes)
arise directly from classical Brownian motion $W$.

\begin{example} \label{E87}
Let $f$ be of class $C^0 (\R)$,   $u \in C^{0,1}(\R_+ \times \R)$.
\begin{enumerate}
\item If $f$ is $C^1$, then  $X=f(W)$ is a $(\shf_t)$-Dirichlet process.
\item $u(t,W_t)$ is an $(\shf_t)$-weak  Dirichlet process,
but not Dirichlet in general.
\item $f(W)$ is not always a Dirichlet process,  
not even  of finite quadratic variation
as shown by  Proposition  \ref{R88}.
\end{enumerate}
\end{example}

The  Example and Remark above easily show that the class
of $(\shf_t)$-Dirichlet processes strictly includes the class
of $(\shf_t)$-semimartingales.

More sophisticated examples  of weak Dirichlet 
processes may be found in the class of the so called
{\it Volterra} type processes, se e.g. \cite{er,er2}

\begin{example} \label{Eer}
 Let $(N_t)_{t\geq 0}$ be an $(\shf_t)$-local martingale,
 $G: \R_+ \times \R_+ \times \Omega  \longrightarrow \R$ 
 a continuous random field
such that $G(t, \cdot)$ is  $(\shf_s)$-{\rm adapted} for each $t$.
Set
$$
X_t = \int_0^t G(t,s) dN_s. $$
Then
 $(X_t)$ is an $(\shf_t)$-weak Dirichlet 
 process with  decomposition $M + A$,
where $ M_t = \displaystyle \int_0^t G(s,s) dN_s.$

Suppose that  $[G(\cdot,s_1);G(\cdot,s_2)]$ exists for any $s_1, s_2$.
With  some additional  technical assumption,
one can show that $A$ is a finite quadratic variation process with
$$ [A]_t = 2 \int_0^t  \left (  \int_0^{s_2}[G(\cdot,s_1);G(\cdot,s_2)]
\circ d M _{s_1}\right ) \circ  d M_{s_2};$$
this iterated  Stratonovich integral can be expressed
as the sum $  C_1(t) + C_2(t) $
where
\begin{eqnarray*}
C_1(t)&=& \int_0^t[G(\cdot,s);G(\cdot,s)]d[M]_s, \\
C_2(t)&= &2\int_0^t \left( \int_0^{s_2}[G(\cdot,s_1);G(\cdot,s_2)]dM_{s_1} \right) dM_{s_2}.
\end{eqnarray*}

\end{example}

\begin{example}
 Take for  $N$ a Brownian motion $W$
and $G(t,s) = B_{(t-s)\vee0}$
where $B$ is a  Brownian motion  independent of $W$.
Then $[A] = \displaystyle \int_0^t (t-s) ds = \frac{t^2}{2}.$
\end{example}

One significant motivation for considering Dirichlet (respectively weak Dirichlet) processes
comes from the study of generalized diffusion processes, typically
solutions of stochastic differential equations with distributional drift.

Such processes were  investigated using 
stochastic calculus via regularization
by \cite{frw1,frw2}. We try to express  here just a guiding idea.
The following  particular case of such  equations  is motivated
by random media modelization:
\begin{equation} \label{EDD}
 dX_t = dW_t + b'(X_t) dt, \quad  X_0 = x_0
\end{equation}
where $b$ is a continuous function.
Typically, $b$ could be the realization  
 of a continuous process, independent of $W$,  
 stopped outside a finite interval.

We shall not recall the precise meaning
 of the solution of (\ref{EDD}).
In \cite{frw1,frw2}  a rigorous
 sense is given to a solution (in the distribution laws) and
existence and uniqueness are established for any initial conditions.

Here we shall just attempt to convince the reader 
that the solution is a Dirichlet process.
For this we define the real function $h$ of class $C^1$ by
$$ h(x) = \int_0^x  e^{-b(y)} dy.$$

We set $\sigma_0 = h' \circ h^{-1}$.
We consider the unique solution in law of the equation
$$ dY_t = \sigma_0 (Y_t) dW_t, Y_0 = h(x_0)$$
which  exists because of classical   Stroock-Varadhan arguments (\cite{sv});
so $Y$ is clearly a semimartingale, thus a Dirichlet process.
The process  $X = h^{-1} (Y)$ is a Dirichlet process since 
 $h^{-1}$ is of class $C^1$.
If $b$ were of class $C^1$, (\ref{EDD})  would be an ordinary
stochastic differential equation, and it could be shown that $X$
is the unique solution of that equation.
In the present case $X$ will still  be the solution of (\ref{EDD}),
considered as a generalized stochastic differential equation.

We now consider the case when the drift is time inhomogeneous as follows:
\begin{equation} \label{EDDI}
 dX_t = dW_t + \partial_x b(t,X_t) dt, X_0 = x_0
\end{equation}
where $b: \R_+ \times \R \rightarrow \R$ is a continuous function of
class $C^1$ in time. Then it is possible to find
a $k: \R_+ \times \R \rightarrow \R$ of class $C^{0,1}$
such that the {\it solution} $(X_t)$ of (\ref{EDDI})
can be expressed as $(k(t,Y_t))$ for some
semimartingale $Y$; so $X$ will 
be an $(\shf_t)$-weak Dirichlet   process.
For this and more general situations, see \cite{rt}.

\subsection{It\^o formula under weak smoothness assumptions} \label{S82}

In this section, we formulate and prove an 
It\^o formula of $C^1$ type.
As for the $C^2$ type It\^o formula, the next
Theorem is stated in the one-dimensional
framework only in spite of its validity in the multidimensional case.

Let ${(S_t)}_{t \ge 0}$ be a  semimartingale and $f \in C^2$.
We recall the classical It\^o formula, as a particular case of  Proposition \ref{P611d}: :
$$
f(S_t) = f(S_0)  + \int_0^t f'(S_s) dS_s + \halb \int_0^t f''(S_s) d[S,S]_s.
$$
Using Proposition  \ref{PITOV3bis} and 
 Definition \ref{DITOV3} (Stratonovich integrals), we obtain
\begin{eqnarray}
  f(S_t)  &=&  f(S_0)  + \int_0^t f'(S_s) dS_s + 
  \halb [f'(S),S]_t \nonumber\\
  \noalign{\vskip-2mm}
&& \label{F83} \\
  \noalign{\vskip-2mm}
 &=& f(S_0) + \int_0^t f'(S) \circ dS. \nonumber
\end{eqnarray}
We observe that in formulae (\ref{F83}), only  the first derivative
of $f$ appears.
Besides, we know that $f(S)$ is a Dirichlet process if $f \in C^1(\R)$.\\
At this point we may ask if formulae (\ref{F83})
remains valid when  $f$ is  in $ C^1 (\R) $ only;
a partial answer will be given in Theorem \ref{T810} below.

\begin{definition} \label{D89}
Let $(S_t)$ be a continuous semimartingale;
set $\hat S_t = S_{T-t}$ for  $t \in [0,T]$. $S$ is called a
{\bf reversible semimartingale} if $(\hat S_t)_{t \in [0,T]}$ 
is again a
semimartingale.\\
\end{definition}

\begin{theorem} \label{T810} (\cite{rv96})
Let $S$ be a reversible semimartingale indexed
by $[0,T]$ and  $f \in C^1(\rz{R})$. Then one has
$$
f(S_t) =f(S_0) +  \int_0^t  f'(S) dS + R_t 
= f(S_0) +  \int_0^t  f'(S) \circ dS
$$
where $R = \frac{1}{2} [ f'(S),S]$.
\end{theorem}

\begin{remark} \label{R810}
After the pioneering work of \cite{by}, which expressed
 the remainder term
$(R_t)$  with the help of generalized integral 
with respect to local time,
 two papers appeared: \cite{fps} in the case of 
 Brownian motion and \cite{fps}
and \cite{rv96} for multidimensional reversible 
semimartingales.  Later, an incredible amount of contributions
on that topic have been published. We cannot give
the precise content of each paper; a non-exhaustive list is 
\cite{bsj, erv,  fp, gp, gr, monu1, monu2}.
Among the $C^1$-type It\^o formulae in the framework 
of generalized Stratonovich integral with respect to
 Lyons-Zheng processes, it is also  important
to quote \cite{lz, lzhe, rvw}.
\end{remark}

\begin{example}
 \label{E811}
\begin{description}
\item{i)} Classical $(\shf_t)$-Brownian motion $W$
is a reversible semimartigale,
see for instance \cite{fps, pa, frw2}.
More precisely
$\hat W_t = W_T + \beta_t + \displaystyle 
\int_0^t \frac{\hat W_s}{T-s} ds$,
where $\beta$ is a $(\shg_t)$-Brownian motion 
and $(\shg_t)$ is the natural
filtration associated with $\hat W_t$.
\item{ii)} Let $(X_t)$ be the solution of the 
stochastic differential equation
$$ dX_t = \sigma(t,X_t) dW_t + b(t,X_t) dt,$$
with $\sigma, b : \R \times \R \rightarrow \R$ 
Lipschitz with at most linear growth,
$\sigma \ge c > 0$. Then $(X_t)$ is a reversible 
semimartingale;
see for instance \cite{frw2}.
Moreover  if $f \in W^{1,2}_{\rm loc} $, it is  
proved in  \cite{frw2} that
$(f(X_t))$ is an  $(\shf_t)$-Dirichlet process.
\end{description}
\end{example}

\begin{proof}

(of  Theorem \ref{T810}).
We use in an essential way
the  Banach-Steinhaus theorem for $F$-spaces; 
see for instance \cite{ds} chap. 2.1.

Define two maps $T_\varepsilon^\pm$
from the $F$-space $C^0(\R)$
 to the $F$-space $\shc ([0,T])$, which consists of 
 all continuous processes indexed by $[0,T]$,
 by
$$
T_\varepsilon^- g = \int_0^{\cdot} g(S_s) 
\frac{S_{s+\varepsilon}-S_s}{\varepsilon}\, ds,
$$
$$
T_\varepsilon^+ g = \int_0^{\cdot} g(S_s) \frac{S_s-S_{s-\varepsilon}}{\varepsilon}\, ds.
$$
These operators are linear and continuous. Moreover, 
for each $g \in C^0$ we have
$$
\lim_{\varepsilon \rightarrow 0} T_\varepsilon^- g 
= \int_0^\cdot g(S) dS,
$$
because of Proposition \ref{PITOV3bis}     which says that
$\int_0^t g(S) d^¯S$ is also an It\^o integral.

Since $\hat S$ is a semimartingale, for the same reasons as above,
\begin{equation} \label{F46}
\int_{T-t}^T g(\hat S) d^- \hat S
\end{equation}
also exists and equals an It\^o integral.

Using Proposition \ref{P24} 3), it follows that
 $\displaystyle\int_0^\cdot g(S) d^+ S$ also exists.

Therefore the Banach-Steinhaus theorem
implies that
$$
g  \mapsto  \int_0^\cdot g(S) d^- S,\quad
g  \mapsto  \int_0^\cdot g(S) d^+ S,
$$
are continuous maps from $C^0(\R)$ to $\shc([0,T])$; by
additivity, so are also
$$
g  \mapsto  [g(S), S], \quad
g  \mapsto  \int_0^\cdot g(S) d^\circ S.
$$

Let $f\in C^1(\R)$, $(\rho_\varepsilon)_{\varepsilon >0}$ be a family
 of mollifiers converging to the Dirac measure at zero.
 We set $f_\varepsilon = f \star \rho_\varepsilon$
where $\star$  denotes convolution. Since $f_\varepsilon$
 is of class $C^2$, by the ``smooth'' It\^{o} formula
 stated at Proposition \ref{P611d} and by Proposition \ref{P24} 1) and 2),
  we have
\begin{eqnarray*}
f_\varepsilon(S_t) &=& f_\varepsilon(S_0) +
\int_0^t f'_\varepsilon(S) dS + \frac{1}{2} [f'_\varepsilon(S),S], \\
f_\varepsilon(S_t) &=& f_\varepsilon(S_0) +
\int_0^t f'_\varepsilon(S) d^\circ S.
\end{eqnarray*}

Since $f'_\varepsilon $ goes to $f'$ in $C^0(\R)$, we can take the limit
term by term and

\begin{eqnarray} \label{F87}
f(S_t) &=& f(S_0) +
\int_0^t f'(S) dS + \frac{1}{2} [f'(S),S], \nonumber \\
\noalign{\vskip-2mm}
&& \\
\noalign{\vskip-2mm}
f(S_t) &=& f(S_0) +
\int_0^t f'(S) d^\circ S.  \nonumber
\end{eqnarray}
Remark  \ref{RITOV4} says that the latter symmetric integral
is in fact a Stratonovich integral.
\end{proof}

\begin{corollary} \label{C814}
If ${(S_t)}_{t\in[0,T]}$ is a reversible semimartingale
and  $g\in C^0(\R)$, then  $[g(S),S]$
exists and
has zero quadratic variation.
\end{corollary}

\begin{proof}

Let $g \in C^0(\R)$ and
let $S = M + V$ be the decomposition of $S$
as a sum of a
 local martingale $M$ and a finite  variation 
 process $V$, such that $V_0 = 0$.
 Let $f \in C^1(\R)$ such that $f' = g$. 
 We know that $f(S)$ is a Dirichlet process
with local martingale part
$$
M_t^f = f(S_0) + \int_0^t g(S) dM.
$$
Let $A^f$  be its zero quadratic variation 
component. Using Thereom \ref{T810}, we have
$$
A_t^f = \int_0^t g(S) dV + \frac{1}{2} [g(S), S].
$$
  $  \int_0^\cdot g(S) dV$ has finite variation, therefore it has
 zero quadratic variation; since so does also $A^f$,
the result follows immediately.
\end{proof}

\begin{proposition} \label{R88}

Let $g \in C^0(\rz{R})$ such that   $g(W)$ is a 
finite quadratic variation process.
Then  $g$ has  bounded variation on compacts.
\end{proposition}
\begin{proof}

Suppose   that $g(W)$ is of finite quadratic variation.
We already know that  $W$ is
 a reversible semimartingale.
By Corollary  \ref{C814}, $[W,g(W)]$ exists and 
it is a zero quadratic variation process.
Since $[W] $ exists, we deduce that
 $(g(W),W)$ has  all its mutual covariations.
In particular $[g(W),W]$ has bounded variation 
because of Remark \ref{R24}.
Let $f$ be such that $f' = g$; Theorem~\ref{T810} implies  that $f(W)$ is a semimartingale.
 A celebrated result of \c{C}inlar,  Jacod,  Protter and  Sharpe
\cite{cjps} asserts that
$f(W)$ is a $(\shf_t)$-semimartingale if and only if $f$
is a difference of two convex functions;
this finally allows to conclude that $g$
has bounded variation on compacts.
\end{proof}

\begin{remark}\label{R815}

Given two processes $X$ and $Y$, the 
covariations $[X]$ and  $[X,Y]$ may exist
even if $Y$ is not of finite quadratic variation.
In particular $(X,Y)$ may not have all its mutual covariations.
For instance, if $X$ has bounded variation, 
and $Y$ is any continuous process,
then $[X,Y] = 0$, see Proposition \ref{P24}~7~b).
A less trivial exemple is provided by
 $X = W$, $ Y = g(W) $ where $g$ is continuous
but not of bounded variation, see
Proposition \ref{R88}.

\end{remark}

\begin{remark}\label{R816} (\cite{fps}).
When $S$ is a Brownian motion,
 Theorem \ref{T810} and Corollary \ref{C814}
are in fact respectively valid  for 
$f \in W^{1,2}_{\rm loc}(\R) $ and $g \in L^{2}_{\rm loc}(\R) $.
\end{remark}

\section{Final remarks} \label{Slast}

We conclude this paper with some considerations about  
calculus related to processes
 having no quadratic variation. On this, the reader can consult
\cite{er2, grv, gnrv}. In \cite{er2} one 
defines a notion of $n-$ covariation  $[X^1, \ldots, X^n]$
of $n$ processes $X^1, \ldots, X^n$ and 
the $n$-variation of a process $X$.

We recall some basic significant results related to those papers.

\begin{enumerate}
\item For a process $X$ having a 
$3$-variation, it is possible to write an It\^o formula of  the type
$$ f(X_t) = f(X_0) + \int_0^t f'(X_s) d^\circ X_s - 
\frac{1}{12} \int_0^t f^{(3)} (X_s) d[X,X,X]_s.$$
Moreover  one-dimensional stochastic differential equations 
driven by
a strong 3-variation were considered in~\cite{er2}.
\item Let $B = B^H $ be a fractional Brownian motion with Hurst index 
$H > \frac{1}{6}$ and $f$ a function of class $C^6$.
It is shown in \cite{grv, gnrv} that
$$f(B_t) = f(B_0) + \int_0^t f'(B) d^\circ B.$$
\item   Using more sophisticated integrals via regularization,
other  types of  It\^o formulae can be
 written  for any $H$  in $]0,1[$;
 see \cite{gnrv}.
\item In \cite{gn}, it is shown that stochastic 
calculus via regularization is
{\it almost pathwise}. Suppose for instance that $X$ is  a  semimartingale
or a fractional Brownian motion,  with Hurst index $H > \frac{1}{2}$;
then its quadratic variation $[X]$ is
a  limit  of $C(\varepsilon, X,X)$ 
not only ucp as in (\ref{SIVR1Cn}),
but also {\it uniformly a.s.}
Similarly, if $X$ is semimartingale and $Y$ 
is a suitable integrand, the  It\^o integral 
$\int_0^\cdot Y dX$ is approximated by
$I^-(\varepsilon , Y, dX)$ not only ucp as in 
(\ref{SIVR1n}), but also uniformly a.s.

\end{enumerate}

\noi{\bf Acknowledgement.}
We wish to thank an anonymous referee and the R\'edac\-tion
of the S\'eminaire 
 for their careful reading of a preliminary
version, which motivated us to improve it considerably.

\def\refname{References}

\bibliographystyle{plain}
\bibliography{Bibliominicours}

\end{document}